\documentclass[a4paper,12pt,onecolumn]{scrartcl} 
\usepackage[a4paper, margin=2.5cm]{geometry}
\usepackage[english]{babel}            %macht deutsche überschriften
\usepackage[latin1]{inputenc}  
\usepackage{amsmath}           %macht
\usepackage{amsfonts}          %       Mathe
\usepackage{amssymb}           %              mächtiger
\usepackage{graphicx}          
\usepackage{caption}
\usepackage{subfigure}
\usepackage{ae}                %macht schöneres ß
\usepackage{typearea}	         %ermöglicht änderung des Seitenspiegels
\usepackage[
  backend=biber,
  style=numeric,
  giveninits=true,   % Vornamen zu Initialen kürzen
  isbn=false,
  maxnames=20,
  maxbibnames=20,
  minbibnames=20
]{biblatex}
\DeclareNameAlias{author}{given-family}
\DeclareNameAlias{editor}{given-family}

\usepackage{lastpage}          %lässt auf die Seienanzahl zugreifen
\usepackage[margin=10pt,font=small,labelfont=bf]{caption} %macht die Bildbeschriftungen richtig
\typearea{16}                    % stellt Seitenspiegel ein
\usepackage{pdfpages}
\usepackage{underscore}
\usepackage{subfigure}
\usepackage{listings}

\RedeclareSectionCommand[
  beforeskip=-2ex,
  afterskip=0.2ex
]{section}

\RedeclareSectionCommand[
  beforeskip=-1.5ex,
  afterskip=0.2ex
]{subsection}
\RedeclareSectionCommand[
  beforeskip=-1.5ex,
  afterskip=0.2ex
]{subsubsection}
\usepackage{booktabs}
\usepackage{setspace}
\usepackage{array}
\usepackage{float}
\usepackage{pdfpages}
\usepackage{lmodern}
\usepackage[backend=biber,style=numeric, isbn=false,maxnames=20,maxbibnames=20,minbibnames=20]{biblatex}
\usepackage{comment}
\usepackage{tabularx} 
\usepackage{scrlayer-scrpage}
\automark{section}

\ihead{}
\ohead{}
\ifoot{\rightmark}              
\ofoot{M. Tesfamarian et al.} 
\usepackage{graphicx} 
\usepackage{graphicx}
\usepackage{placeins}
\usepackage{hyperref}
\usepackage[T1]{fontenc}

\makeatletter
\let\@fnsymbol\@arabic
\makeatother

\title{\Large \textbf{Modelling Human Skin Morphology and 
Simulating Transdermal Transport of 50 Chemicals
}}

\author{\normalsize  Milana Tesfamarian\footnotemark[1], 
\normalsize  Arne N{\"a}gel\footnotemark[2],
\normalsize  Michael Heisig\footnotemark[2],
\normalsize  Gabriel Wittum\footnotemark[1]
}
\date{} % kein Datum

\begin{document}
\pagenumbering{arabic}
\maketitle
% Fußnotentexte unter dem Titel
\footnotetext[1]{King Abdullah University of Science and Technology (KAUST), Computer, Electrical and Mathematical Sciences and Engineering (CEMSE) Division, Thuwal 23955-6900, Saudi Arabia.\\ Correspondence: milana.tesfamarian@kaust.edu.sa}
\footnotetext[2]{Goethe University, Simulation and Modelling, Kettenhofweg 139,  Frankfurt, Germany.}

\begin{abstract}
\small                % kleinere Schrift
\setstretch{1.5}
SUMMARY: The goal of enhancing effective and controlled delivery of active ingredients, focusing on developing robust computational simulations, is gaining importance. The development of mathematical models, combined with an understanding of a system's underlying physics and chemistry, is a powerful tool for gaining more precise and quantitative insights into the delivery of active ingredients over time and space. To establish the transdermal patch for practical use, we must understand the chemical transport mechanisms from the patch to the skin. Furthermore, it is essential to consider different body regions with various parameters in the sublayers. In addressing this issue, a mathematical model is formulated that incorporates key parameters, including skin characteristics and the physicochemical properties of the transdermal patch. We investigated the diffusion and permeation behavior of 50 chemicals and peformed numerical simulations using parameters that represent the chemical absorption and diffusivity within the specific sublayers.

\setlength{\parskip}{2em}


\end{abstract}            

\newpage

\section{Introduction}\label{subsec1}
The human skin is the body's largest organ. It is an essential component of the human body and has a complex structure that consists of many individual layers. The skin's primary function is to protect internal organs and muscles from pathogens and other environmental influences \cite{yousef2017anatomy}. The protective mechanisms are realized through the layer hierarchy of the skin, particularly by the three main skin layers:
\begin{enumerate}
\item Epidermis 
\item Dermis 
\item Subcutis 
\end{enumerate} 
In these three layers, essential events and processes occur, from blood supply to oxygen exchange. The outermost layer is the stratum corneum, which is the main barrier to aqueous substances and protects the body from dehydration. The epidermis consists of sublayers, which can be collectively defined as the viable epidermis, excluding the stratum corneum (SC). The dermis, which lies beneath the epidermis, contains hair follicles, sweat glands, and blood vessels. The dermis is also divided into individual layers:
\begin{enumerate}
\item Papillary dermis
\item Reticular dermis
\end{enumerate}
The interface between the epidermis and the dermis is called the dermoepidermal junction, and the wave-shaped layer beneath the epidermis is called dermal papillae \cite{b}. These finger-shaped structures contribute significantly to the mass exchange of oxygen and nutrients \cite{baqar2019dermal}. Furthermore, the dermal papillae contain capillary loops and sensory nerve endings. The deepest layer is called the subcutis. It consists of connective tissue that maintains elasticity and ensures the skin's resistance to tearing. The following picture shows a 2D microscopic image with relevant layers, combined with an abstract model. This includes the stratum corneum (SC), viable epidermis (VE), papillary dermis (PD), and reticular dermis (RD). The abstract model is illustrated as a layer-by-layer stack, like a \glqq sandwich\grqq{}. The epidermis and the upper part of the dermis, more specifically the dermal papillae, are the primary topics of discussion in this paper.
\begin{figure}[htbp]                                 
\centering
\includegraphics[scale=0.6]{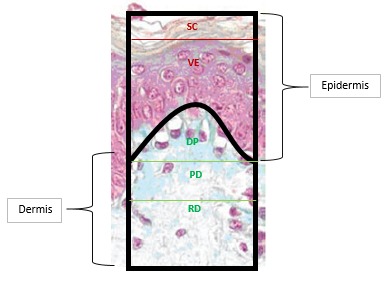} 
\caption[]\\{Microscopic image \cite{ga} and abstract schematic representation to demonstrate the structural organization of the different skin layers} 
     \label{fig:genisiert}
     \end{figure}
\section{Methods and Materials}
\subsection{Comparison of Scientific Articles about the Skin}
A realistic skin model is developed based on microscopic data by comparing the measured data and test procedures. All measurement techniques yield different mean values and ranges. Nevertheless, these techniques are useful and provide fast orientation values. Several aspects are responsible for the differences in the measurements. One of the factors that contributes to the variability of the data is the age gap that exists between subjects. For example, \cite{n7} investigated subjects aged between 19 and 30 years. In contrast, \cite{n4} examines the density of dermal papillae per unit area. The number of dermal papillae was significantly lower in older subjects than in younger ones. Additionally, the thickness of their skin was investigated during the winter. According to measurement data in \cite{n1}, skin thickness varied slightly between female and male subjects at an average age of 77.5 years. \cite{n77} measured the skin thickness of Korean men and women. Comparing the measured epidermal thicknesses between men and women reveals that, in numerous skin locations, male epidermis exhibits greater thickness. Environmental conditions also influence the skin thickness. People living in sunnier regions tend to have thicker skin than those who live in less sunny areas \cite{g}. The study included subjects with African or Caucasian ancestry who live in Europe. In African chest skin, the dermal-epidermal interface was three times longer than in Caucasian chest skin. Furthermore, differences were observed not only in skin thickness but also in the cellular interaction between the epidermis and dermis. According to \cite{n3}, there is a larger dermal epidermal junction when the skin is exposed to high levels of sunlight.
\begin{figure}[h]
\centering
\includegraphics[width=0.45\linewidth]{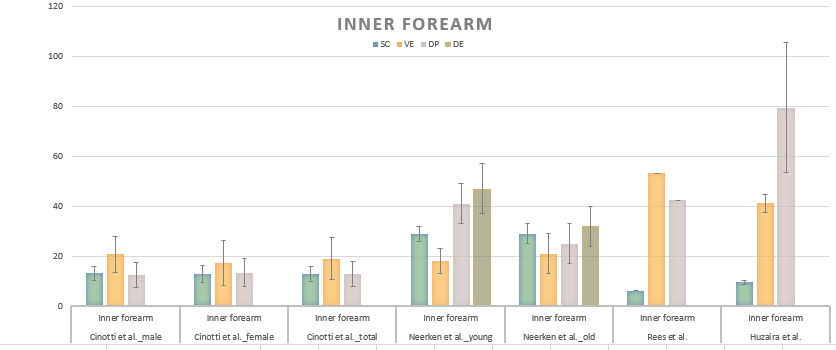}
\vspace{1em} 
\includegraphics[width=0.45\linewidth]{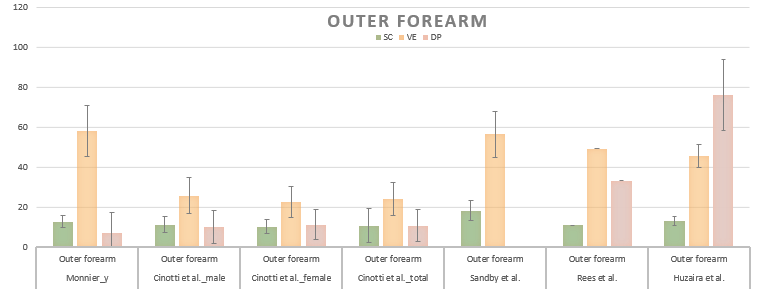}
\vspace{1em} 
\includegraphics[width=0.35\linewidth]{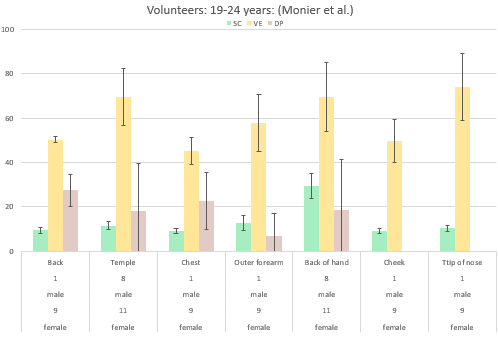}
\includegraphics[width=0.35\linewidth]{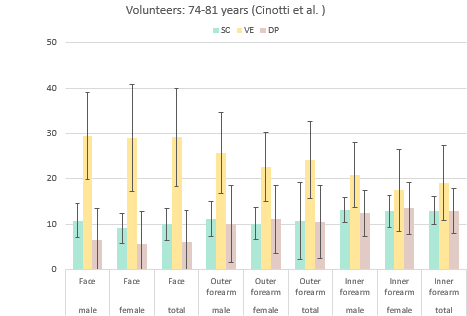}
\includegraphics[width=0.35\linewidth]{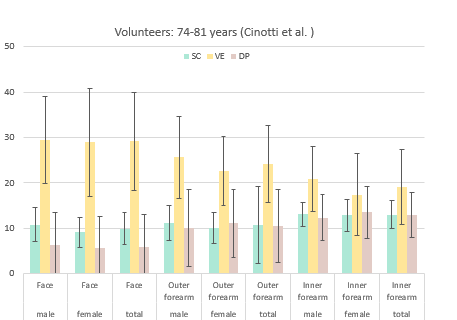}
\includegraphics[width=0.25\linewidth]{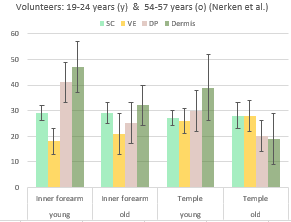}

\caption{Overview of microscopic data and comparison of structural characteristics across various skin regions}
\label{fig:general}
\end{figure}

In summary, the following attributes influence the skin thickness: 
\begin{itemize}
\item Ethnicity
\item Age
\item Gender
\item Health
\item Environment, Climate
\end{itemize} 

\begin{figure} [H]
\begin{center}
\fbox{\includegraphics[scale=0.55,angle=270]{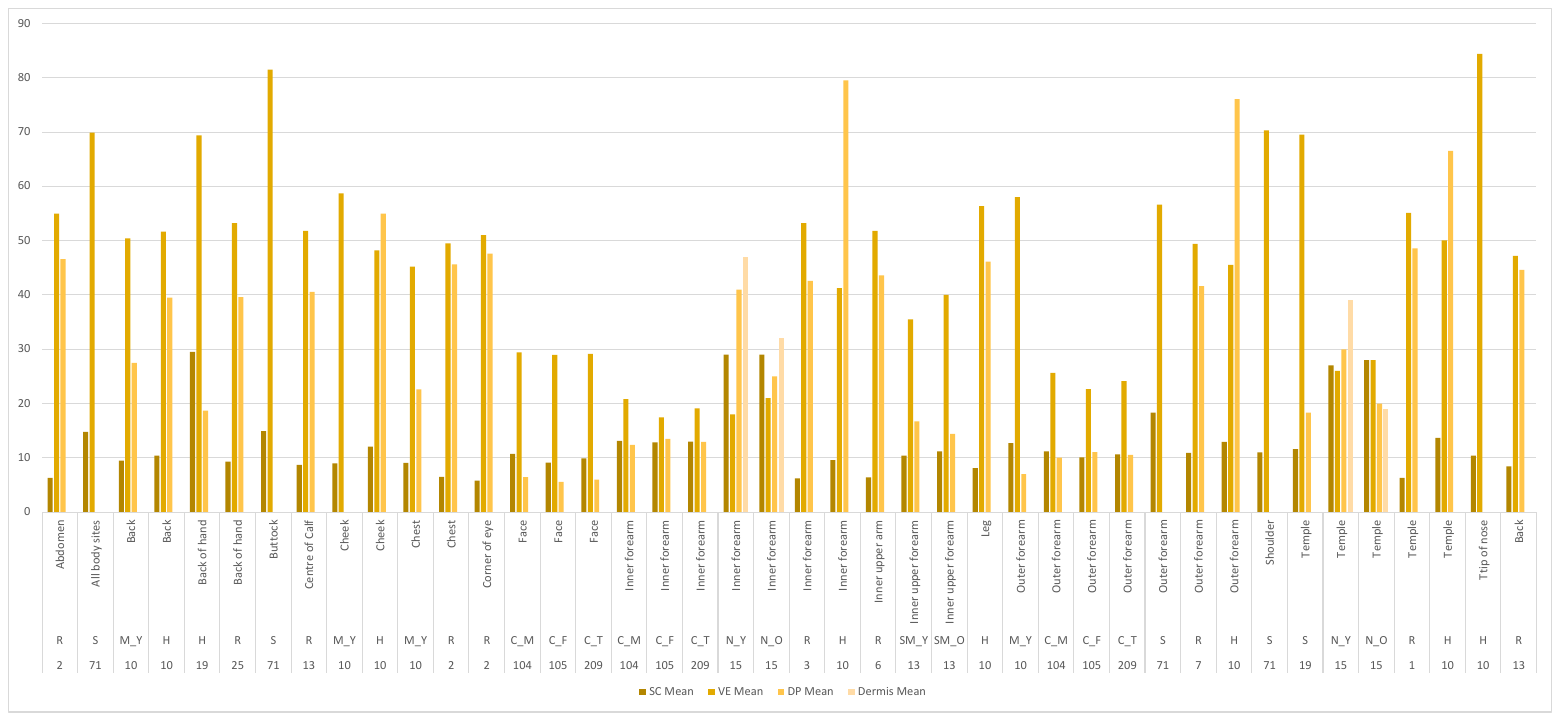}} %% adjust angle and scale appropriately
\caption{
The bar chart includes mean values for different skin layers, with skin thickness measurements [in $\mu$m] of various body regions from $n$ study participants from diverse scientific articles, enabling comparison across anatomical sites. 
}
\label{fig:TEg}
\end{center}
\end{figure}

\noindent
The bar charts presented in Figure \ref{fig:general} and Figure \ref{fig:TEg} provide an overview of the microscopic skin thicknesses of various body regions, obtained from different scientific publications.

\subsection{Finite Dosing}
\label{subsec:masse}
In dermal drug delivery studies, finite dosing experiments are frequently conducted. This concept is commonly used to simulate real-world scenarios that involve the application of creams to the skin. We demonstrate the finite dosing case in which the deposition layer contains one chemical. The deposition layer, which is positioned above the SC,  has a thickness of 50~\textmu m, corresponding to the thickness of the transdermal patch used in \cite{Tfentanyl}. Over time, the chemical in the deposition layer decreases until it is completely depleted. 

\subsection{Model Geometry}
\begin{figure}[htbp]                                  
\centering
\includegraphics[scale=0.5]{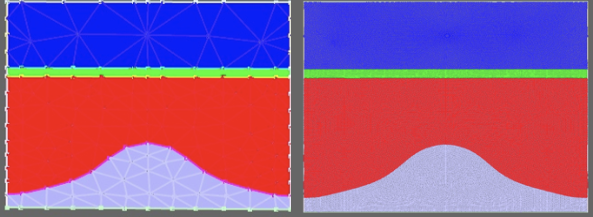} 
\caption{Mesh- older skin - skin region chest - finite dose - 2D}
\label{fig:ref5depos}
\end{figure}

Realistic and computable skin meshes in 2D and 3D are developed based on scientific publications by comparing and analyzing measurement techniques, microscopic data and images \cite{n6}\cite{n4}. Using these data, we developed functions describing the three-dimensional surface of the dermal papillae. The following functions are:
\begin{eqnarray}
f(x)=\frac{h}{2}\cdot sin(2\pi \cdot \frac{ x}{a}-\frac{1}{2}\pi)+\frac{h}{2}
\label{paa1}
\end{eqnarray}
\begin{eqnarray}
g(y)=\frac{\frac{h}{2}\cdot sin(2\pi \cdot \frac{y}{a}-\frac{1}{2}\pi)+\frac{h}{2}}{h}
\label{paa2}
\end{eqnarray}
\noindent
Here, \textit{h}  represents the thickness of the dermal papillae and \textit{a} the period, with $a\mathrm{=}200$ \textmu m (old skin) or $a\mathrm{=}150$ \textmu m (young skin). Figure \ref{fig:ref5depos} displays the initial coarse skin mesh (left), consisting of 171 vertices, 460 edges, and 290 triangles, and the refined mesh (right) after five refinement steps. The blue layer indicates the DEPOS layer, while the green layer illustrates the SC. The red layer represents viable epidermis, and the sky-blue layer the dermis.

\begin{figure}[htbp]                                 
\centering
\includegraphics[scale=0.6]{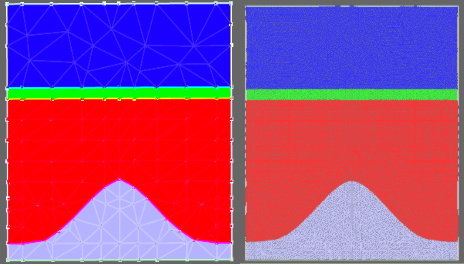} 
\caption{Mesh- young skin - skin region chest - finite dose - 2D}
\label{fig:ref6depos}
\end{figure}

\subsection{Model Equation}
\noindent
The transport of chemicals can be described using Fick's second law, as this is a transient diffusion process in which the dependence of the concentration $c$ on location $x$ and time $t$ is considered \cite{mdiff}. If a layered model exists, the solubility of the substance in the respective layers leads to variations in concentration. Through Nernst's distribution law in \eqref{eq:neun}, the solubility of a substance in two adjacent phases can be characterized. Moreover, the partition coefficient $K$ can be described as the ratio of concentrations in the respective layers $i$ and $j$. However, this is not yet the final equation for the partition coefficient. Since \textit{K} is not always accessible, the reference concentration \textit{u}, which is known as the chemical potential, can also be used.
\begin{equation}
\begin{split}
K_{ij}=\frac{c_{i}}{c_{j}}
\end{split}
\label{eq:neun}
\end{equation}
\begin{equation}
\begin{split}
K_{ij}= \frac{c_{i}}{c_{j}}=\frac{K_{i}u}{K_{j}u}
\end{split}
\end{equation}
\noindent
In conclusion, the following general model results in relation for a substance, so that the index \textit{i} or \textit{j}  can be omitted:
\begin{eqnarray}
\nabla \cdot (D\nabla Ku) &=& \frac{\partial(K(x)\cdot u)}{\partial t} \label{true2}
\end{eqnarray}

\subsection{Conservation of Mass}
Another key aspect is to ensure that mass conservation is maintained in the discrete case. The conservation law underlying the following equations is satisfied at the discrete level. It must be integrated over the interface to calculate the masses in the respective skin layers. Using the formula in \eqref{eq:eins}, the trapezoidal rule, the released mass is calculated.
 \begin{equation}
 \begin{split}
m(t_{k}):=\int_{0}^{t_{k}} f(t)\,dt \approx  \sum_{i=1}^{k} \frac{f(t_{i}) + f(t_{i-1})}{2} \cdot (t_{i}-t_{i-1}) 
 \end{split}
 \label{eq:eins}
 \end{equation}
In this case, $f(t)$ is the mass flow. This is calculated with the following formula:
\begin{equation}
\begin{split}
f(t)=\frac{1}{|\tau|}\int_{\tau} -KD \nabla u(t,x) \cdot  \vec{n}~dA,    \text{ where }  \tau =  (t_{i}-t_{i-1}) 
\end{split}
 \label{eq:zwei}
\end{equation}
\noindent
To check whether the physical principle of mass conservation is given, the total mass is calculated using the formula $M_{total}= M_{SC} + M_{VE} + M_{DE} + M_{releasedmass}$, while $M_{SC}$, $M_{VE}$ and $M_{DE}$ are the calculated masses per area in the respective skin layers.

\subsection{The Finite Volume Method for Solving the Diffusion Equation}
\label{sec:finit}
\begin{figure}[htbp]                                
\begin{center}                                       
\includegraphics[width=0.4\linewidth ]{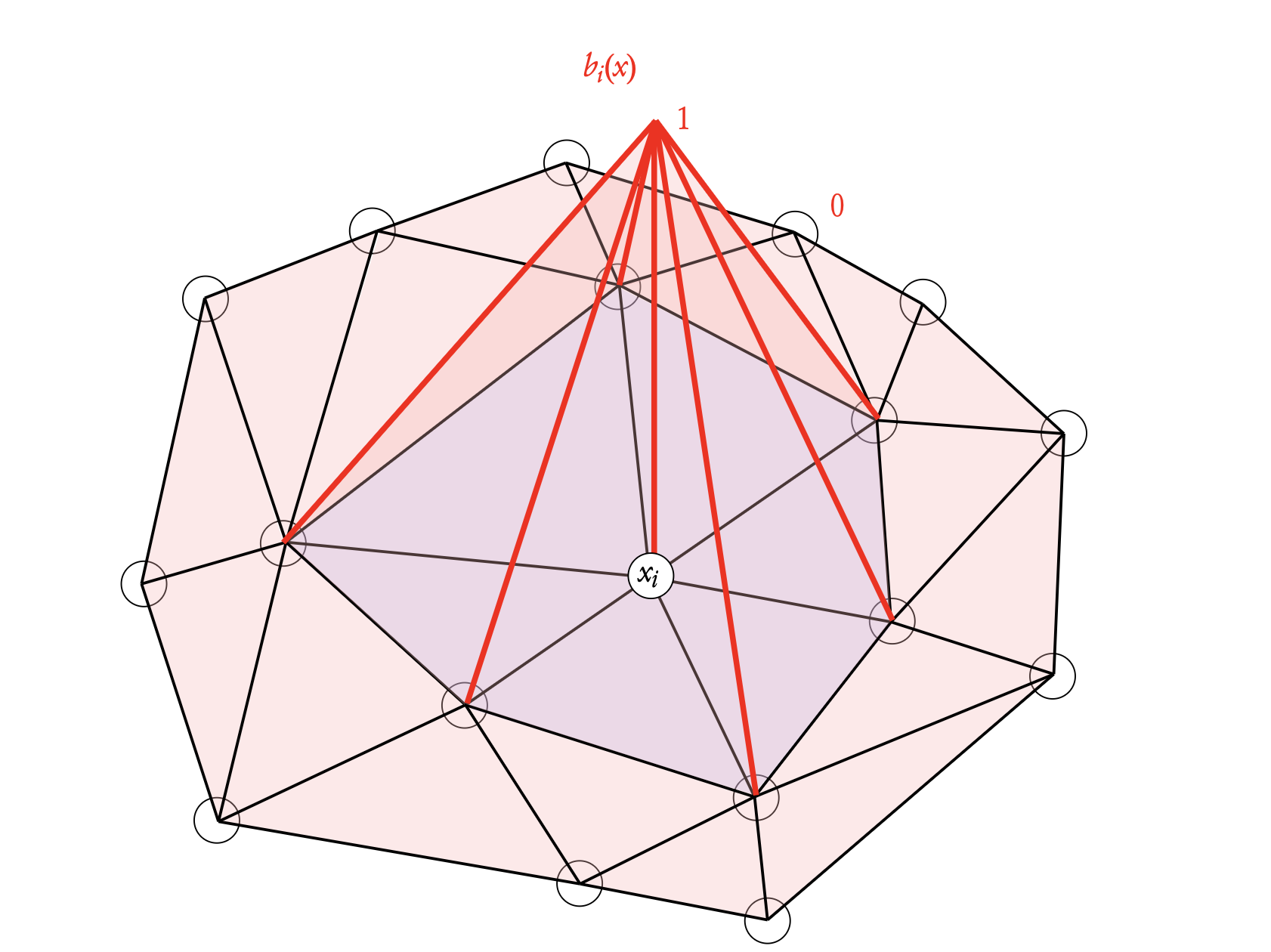}       
\caption{A piecewise linear basis function on a triangular grid }     
\label{fig:vor}              
\end{center}
\end{figure}

To compute the equations numerically, they must be approximated on the grid. For this purpose, the second-order of the finite volume method is described below \cite{knabner2013numerik}. For a given volume $V$, the equation is expressed as follows: 
\begin{eqnarray}	
\int_V \frac {\partial u}{\partial t} = \int_{\partial V}  D\nabla u \cdot \vec {n} \label{diff}
 \label{eq:partdrei}
\end{eqnarray}	
This holds for an arbitrary volume $V$. To construct appropriate integration volumes, the domain $\Omega$ is assumed to be covered by a triangulation $T$. At each vertex, the midpoints of all incoming edges and the barycenters of the triangles are marked. The edge midpoints are then connected to the barycenters as illustrated in Figure \ref{fig:vor}. This process creates a box around each grid point. For each of these boxes, the equation \eqref{eq:partdrei} is fulfilled.

The Gauss theorem is now applied to the left-hand side of the equation. The result is as follows:
\begin{eqnarray}
\begin{split}
\int_{\partial{B_{i}}} (DK\nabla u)\cdot \vec{n}\,dx= \int_{{B_{i}}} \frac{\partial( Ku)}{\partial t}\,dx\\ \label{gh44}
\end{split}
\end{eqnarray}
Therefore, a basis representation of $u(x)$ is created:
\begin{eqnarray}	
u(x,t) = \sum\nolimits_{j}  u_{j}(t) \chi_{j}(x).
\end{eqnarray}
\noindent
With linear finite elements, the $\chi_{j}(x)$ are so-called \textit{hat functions}. A hat function is a piecewise linear function defined as:
\begin{eqnarray}	
\begin{split}
\chi_i(x) =
\begin{cases}
\frac{x - x_{i-1}}{x_i - x_{i-1}}, & x_{i-1} \leq x \leq x_i, \\
\frac{x_{i+1} - x}{x_{i+1} - x_i}, & x_i \leq x \leq x_{i+1}, \\
0, & \text{otherwise}.
\end{cases}
\end{split}
\end{eqnarray}
\noindent
\newline
This finally results in a linear system of equations that can be calculated in the following form:
\begin{align}	
\sum \nolimits_{j}\int_{\partial{B_{i}}}DK \nabla \chi_{j}(x) \vec{n} dx u_{j}(t)= \sum \nolimits_{j} \int_{Bi} K \chi_{j}(x)dx\frac{\partial u_{j}(t)}{\partial t}  \label{gh45}
\end{align}
Applying an implicit Euler discretization to the time derivative, we take into account that the basis function $\chi_i(x)$ can be calculated. This yields to (\ref{equation_MVec}) with a time step $\tau >0$, the stiffness matrix $A$, and the mass matrix $M$. The left-hand side of the equation can be summarized as $(M+\tau \cdot A) \cdot \vec{u}(t+\tau) $, and we denote the right-hand side by $M\vec{u}(t)$. 
The system of equations takes the following final form after applying the Euler method.
\begin{eqnarray}	
\begin{split}
(M+\tau\cdot A)\cdot \vec{u}(t+\tau) &=& M\vec{u}(t) 
\label{equation_MVec}
\end{split}
\end{eqnarray}
Finally, (\ref{equation_MVec}) can be solved using the multigrid method. The LIMEX method (\textit{Linearly Implicit Method Discretization combined with Extrapolation}) is used for adaptive time stepping \cite{ux}\cite{uu}. If the time steps are too large, the algorithm reduces the step sizes. The following differential equation with initial value $u_{0}$ is given:
\begin{eqnarray}
u^{\prime}&=&f(u)
\end{eqnarray}
The boundary conditions are essential to describe the diffusion behavior at the boundaries $\partial \Omega$. The simulations performed in this study employ the Dirichlet boundary condition. The chemical is transported into the underlying skin layers until it reaches the lowest interface $u_{BOT}$~$\mathrm{=}$~0, where a sink has been implemented.

\subsection{Numerical Implementation}
For solving large and sparse systems of linear equations in linear time, we work with the multigrid method \cite{hackbusch}. By using iterative methods such as the Gauss-Seidel method or the Incomplete LU Decomposition, we are able to minimize time complexity and to improve convergence efficiency. The finer the grid, the better the precision of the simulation results. This indicates that increased refinement can maintain mass conservation, but more computing time will be required. For this reason, it is advantageous to apply numerical methods in which the level of refinement can be selected. 

\begin{figure}[htbp]                                 
\centering
\includegraphics[scale=0.55]{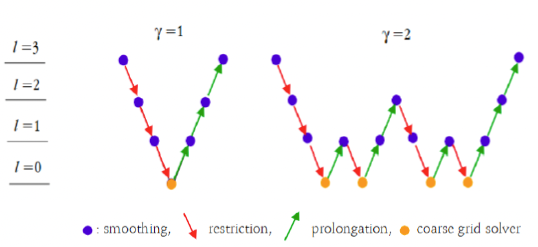}
\caption{Multi-grid cycles. Left: V-cycle (y=1), right: W-cycle (y=2)}   
\label{fig:pic_multigri}
\end{figure}

\subsection{Software}
\noindent
The simulations were performed by programming in UG4 \cite{ug4}. UG4 is a simulation framework written in C++. It is a cross-platform tool and supports Linux, macOS, and Microsoft Windows. UG4 can solve differential equations using unstructured hybrid grids in 1D, 2D, and 3D, and allows fast and efficient numerical solution methods. Furthermore, UG4 is supported by several powerful tools, such as the meshing software ProMesh \cite{proM}. These tools can be operated using the Lua scripting language. In addition, simulation results can be saved and interpreted by applying widely used visualization toolkits.

\subsection{Model Parameters}
\begin{figure}[htbp]
\centering   
\includegraphics[scale=0.6]{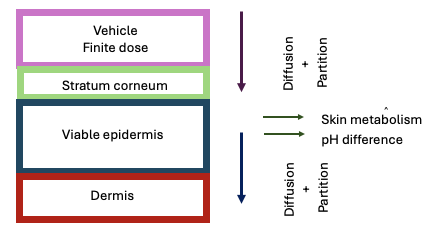}                         
\caption{A schematic view of the skin illustrating the various layers used in the simulation.}   
\label{fig:skizze}
\end{figure}

To compute the diffusion in the sublayers and to consider the interfaces between the layers, it is necessary to define $K$ and $D$ (Figure \ref{fig:skizze}). 

Using formulas from research articles \cite{ellison}\cite{dancik}\cite{selzer}, missing parameters were defined and calculated. For example, the diffusion coefficient of the viable epidermis was calculated according to the formula \cite{dancik}:
\begin{equation} 
 D_{free}=10^{(-4.15-0.6555\cdot logMW)}
 \end{equation}
 $D_{SC}$ was calculated for the chemical thioglycolic acid using the SC height and lag time by applying the formula  (\ref{eq:SC_equation}) \cite{kruczek2005analytical}. Despite the lack of direct data on thioglycolic acid, experimental investigations show a range of lag times and that thioglycolic acid penetrates the skin quickly. According to these data, the lag time for thioglycolic acid can be between one and five hours \cite{Thio}. However, a specific experimental measurement would be necessary for a more accurate determination. 
  \begin{equation}
D_{\text{SC}} = \frac{H_{\text{SC}}^{2}}{6 \cdot t_{\text{lag}}}
\label{eq:SC_equation}
 \end{equation}

\subsection{Chemicals}
After simulating the diffusion of chemicals, certain chemicals were categorized and observed. The selection was made based on their specific properties and characteristics, such as molecular weights, pH values, coefficients, and simulation findings. Although some of the chemicals here can have toxic or harmful effects on human health, they are still used in commercial products at low concentrations. The chemical triclosan has a molecular weight of 289.542 Da. It has antibacterial properties and is therefore used in products such as toothpaste, mouthwash, hand sanitizer, and surgical soaps. In September 2016, triclosan was banned in cosmetic products due to its high absorption through human skin and oral mucosa \cite{triclosan}. The chemical nitrobenzene has a molecular weight of approximately 123.109 Da. In the mid-19th century, this chemical was used to perfume inexpensive curd soaps, replacing the more expensive almond oil \cite{nitro}. Nitrobenzene was later banned from cosmetic products as a result of its toxic effects. Benzylideneacetone (MW: 146.186 Da) is used as an additive preservative in the food and cosmetic industry and is classified as less toxic \cite{acetone}. Basic Red 76 (MW: 371.8 Da) is a chemical used in cosmetic products and as an ingredient in hair coloring \cite{cherian}. It is considered safe and approved by regulatory bodies such as the FDA and the EU. However, it can cause skin irritation or allergic reactions in sensitive individuals. Propylparaben (MW: 180.2 Da) is commonly found in various cosmetic products, including shampoos, lotions, and creams, to prevent them from spoiling \cite{propyl2}. However, distinct types of parabens can cause serious health issues such as allergies and hormone disruption \cite{propyl}.  
\noindent
2-Ethylhexyl acrylate (MW: 184.27535 Da) is a widely utilised acrylate monomer in the synthesis of polymer resins, serving various applications including adhesives, latex paints, cross-linking agents, textile and leather finishes, as well as paper coatings \cite{Ethy2}. In laboratory trials, the administration of 2-ethylhexyl acrylate resulted in dermal alterations in mice, including irritation, scaling, scabbing, epidermal thickening, and heightened cellular proliferation \cite{Ethy21}. Naphthalene (MW: 128.1705 Da) was first registered as a pesticide in the United States in 1948. It is used in the production of dyes, resins, and solvents, and also has applications in pest control. However, naphthalene can cause health issues, including respiratory problems and, in severe cases, damage to red blood cells \cite{Naph0}.  Due to its antiseptic and disinfectant effects, Resorcinol (MW: 110.1 Da) is used in various skin treatments for conditions such as acne, seborrheic dermatitis, and eczema \cite{Rc}.

\section{Results and Discussion}
Only the most distinct behaviors and noteworthy patterns are presented here. The meshes (in 2D and 3D with DEPOS layer) were refined five times. With finite dosing, the mass of the individual skin layers decreases to zero over time, as the top layer diminishes and the chemical penetrates the inner skin layers. The total mass remains constant. The released mass indicates how much of the chemical has been absorbed by the skin. The following sections contain the description and plots of certain simulation results (\ref{subsec:Simresults}).
On the x-axis, the time is displayed in days, and the mass per area is indicated in \mbox{\textmu g\slash\textmu m\textsuperscript{2}} on the y-axis. For each time point, the masses of the DEPOS layer and the three skin layers (SC, VE, and DE) are plotted, along with the released and total mass. In addition, the maximum mass per area (M\textsubscript{max}) and the associated time (T\textsubscript{max}) are evaluated. Furthermore, many experimental tests are performed using human skin samples with a duration of up to 200 hours \cite{Neil}. It appeared ineffective and impractical to look at skin penetration after 20 days, because skin samples often cannot withstand the chemicals. Thus, samples were analyzed between 2 and 16 days.   

\subsection{Comparison between Young and Old Skin}
\begin{figure}[htbp]                                 
\centering
\includegraphics[scale=0.2]{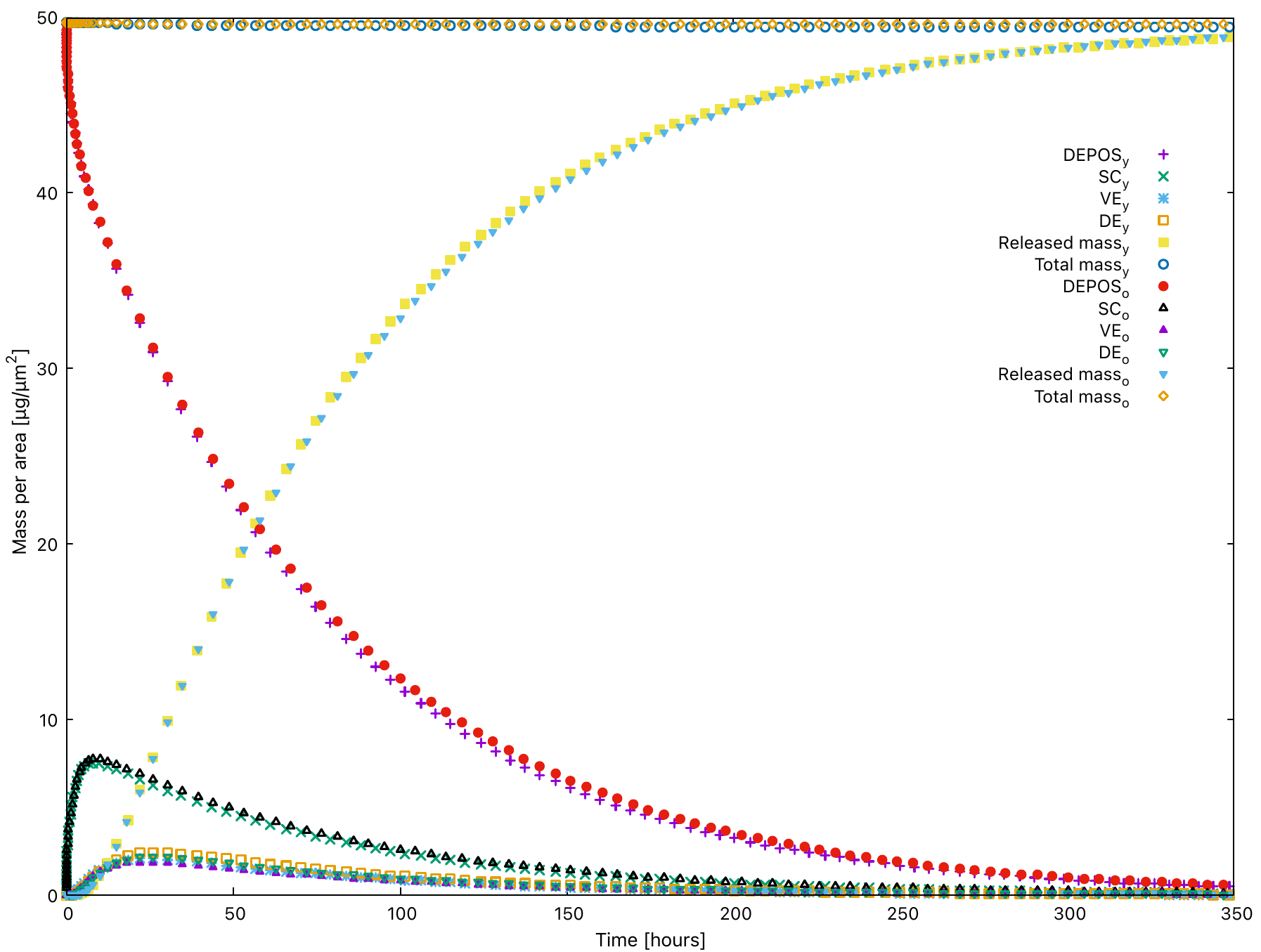}
\caption{ Plots of the simulation results for benzyl bromide using young and old skin meshes in 2D}   
\label{fig:OLDYOUNG}
\end{figure}

According to \cite{OY}, we distinguished between young and old skin by changing the length and width of the dermal papillae. Younger skin has longer dermal papillae (Figure \ref{fig:ref6depos}) while older skin has flatter dermal papillae (Figure \ref{fig:ref5depos}) \cite{blair2020skin}. Despite the different geometric representations for older and younger skin, the simulation results show minimal differences (Figure \ref{fig:OLDYOUNG}). However, the shape of the dermal papillae alone cannot explain the differences between young and old skin. Other factors, such as the thickness and shape of the sublayers, could also play a significant role. With advancing age, the skin becomes thinner and the dermal capillary density decreases \cite{vybohova2012quantitative}. Additionally, a reduction in dermal capillary density could affect transdermal drug permeability.

\subsection{The Comparison of Chemicals}
\subsubsection{Basic Red 76 and 2-Ethylhexyl acrylate}
\begin{figure}[htbp]                                 
\centering
\includegraphics[scale=0.22]{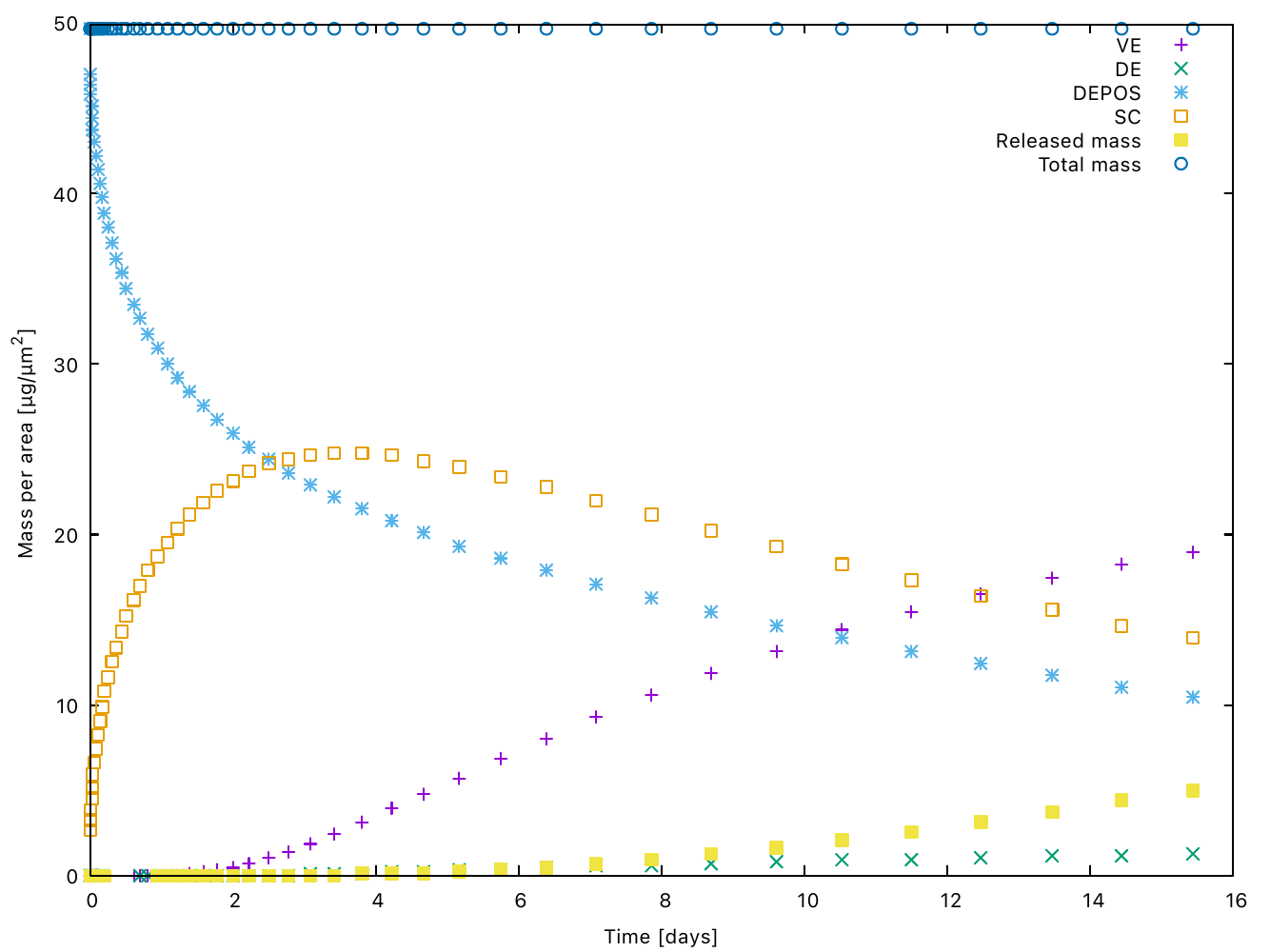}%{plot_2_complete.pdf}
\includegraphics[scale=0.22]{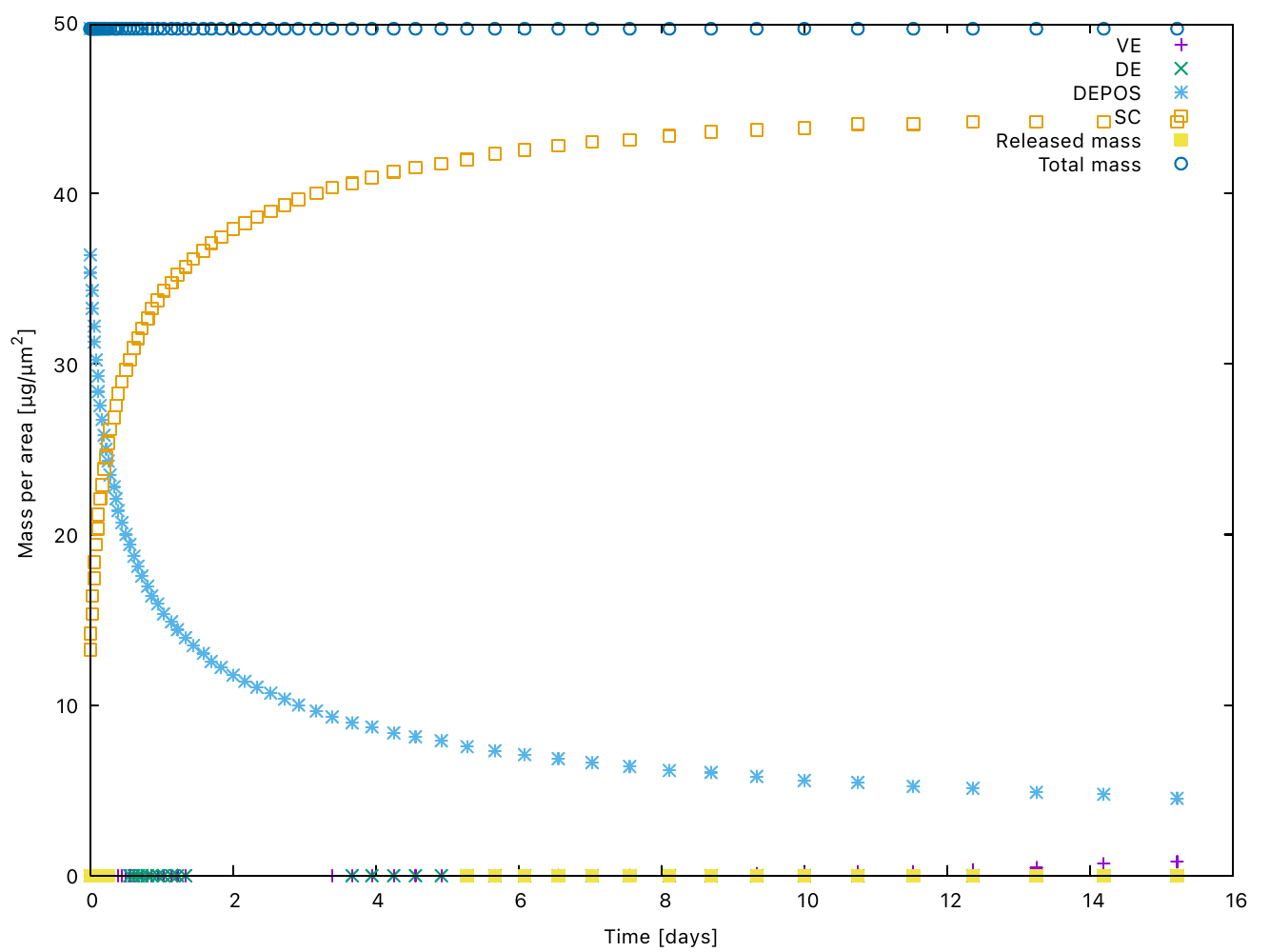}
\caption{Plots of the simulation results for basic red 76 (left) and 2-ethylhexyl acrylate (right) in the abdominal skin mesh} 
\label{fig:B2}
\end{figure}

At $t = 0$, the simulation represents a fully saturated DEPOS layer containing one specific chemical. After $t>0$, the chemical diffuses into the underlying SC layer. The SC absorbs the chemical basic red 76 with a maximum mass per area of 24.44 \textmu g\slash \textmu m\textsuperscript{2} in approximately 3.5 days (Figure \ref{fig:B2}). This equals approximately half of the total mass in the DEPOS layer. The VE absorbs the maximum mass per area for basic red 76 in approximately 25 days. Furthermore, the maximum mass in the VE is 23.54 \textmu g\slash \textmu m\textsuperscript{2}, while in the dermis it is 1.639 \textmu g\slash \textmu m\textsuperscript{2}. In contrast, 2-ethylhexyl acrylate penetrates the SC more effectively than basic red 76. The enhanced absorption is due to the higher partition coefficient of 2-ethylhexyl acrylate in the SC ($K_{SC}=170$) in comparison to that of basic red 76  ($K_{SC}=18$). The SC absorbs a maximum mass per area of \mbox{44.28 \textmu g\slash \textmu m\textsuperscript{2}} in approximately 13 days. 2-Ethylhexyl acrylatepenetrates almost entirely from the DEPOS layer to the SC. The VE absorbs approximately half of the mass of the DEPOS layer, which is around 26 \textmu g\slash\textmu m\textsuperscript{2}. The dermis absorbs 2-ethylhexyl acrylate with a maximum mass per area of 2.28 \textmu g\slash\textmu m\textsuperscript{2}. In the mesh skin, no blood vessels were included within the dermal papillae. Here, the scientific article \cite{kretsos} should be referenced, in which the capillaries were developed as a microscopic model. Moreover, we observe minimal diffusion in the dermis, which reflects a realistic case, since the capillary loops in the dermal papillae absorb most of the chemicals.

\subsubsection{Naphthalene and Propylparaben}
\begin{figure}[htbp]                                        
\centering
\includegraphics[scale=0.22]{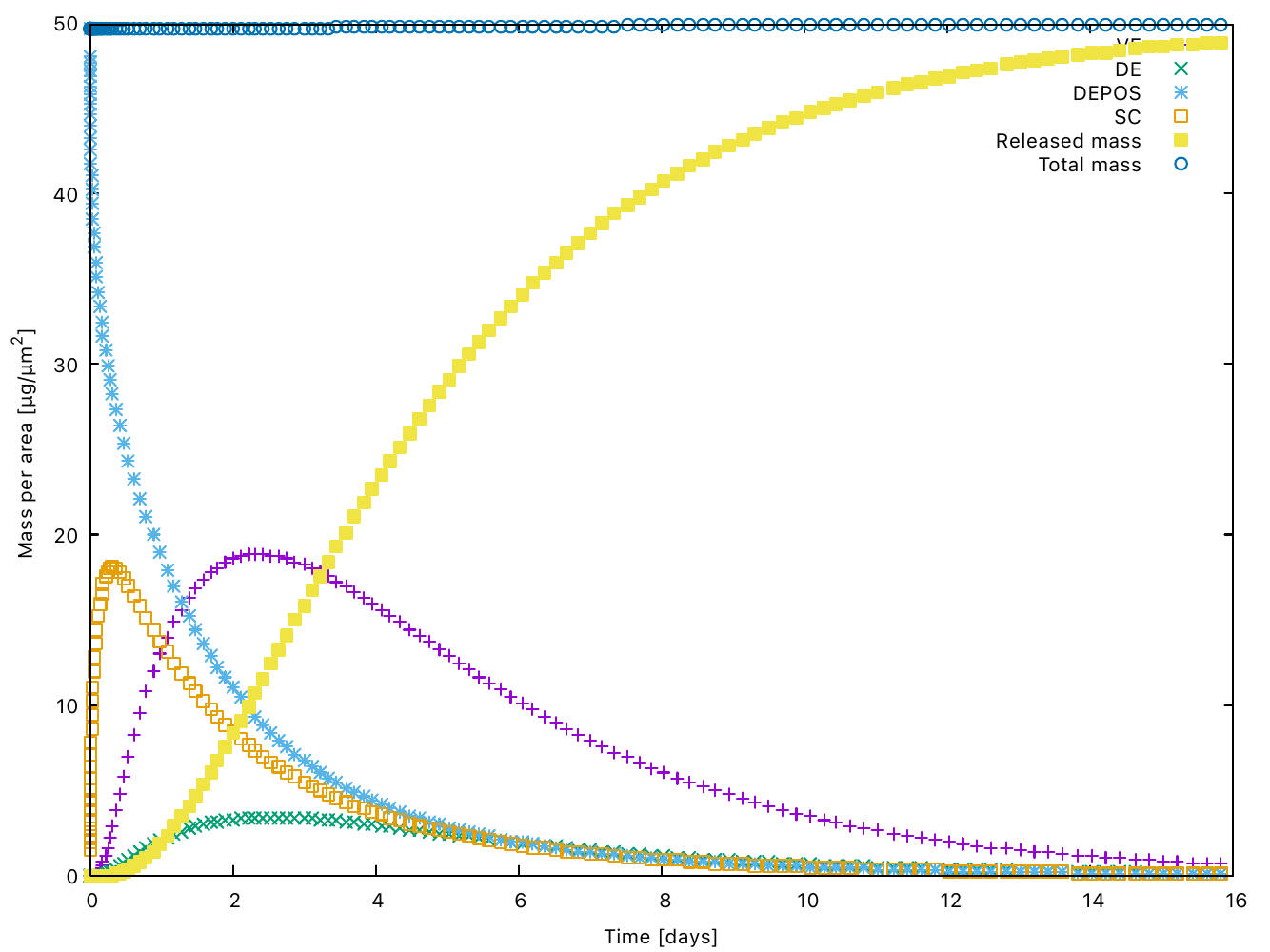}
\includegraphics[scale=0.22]{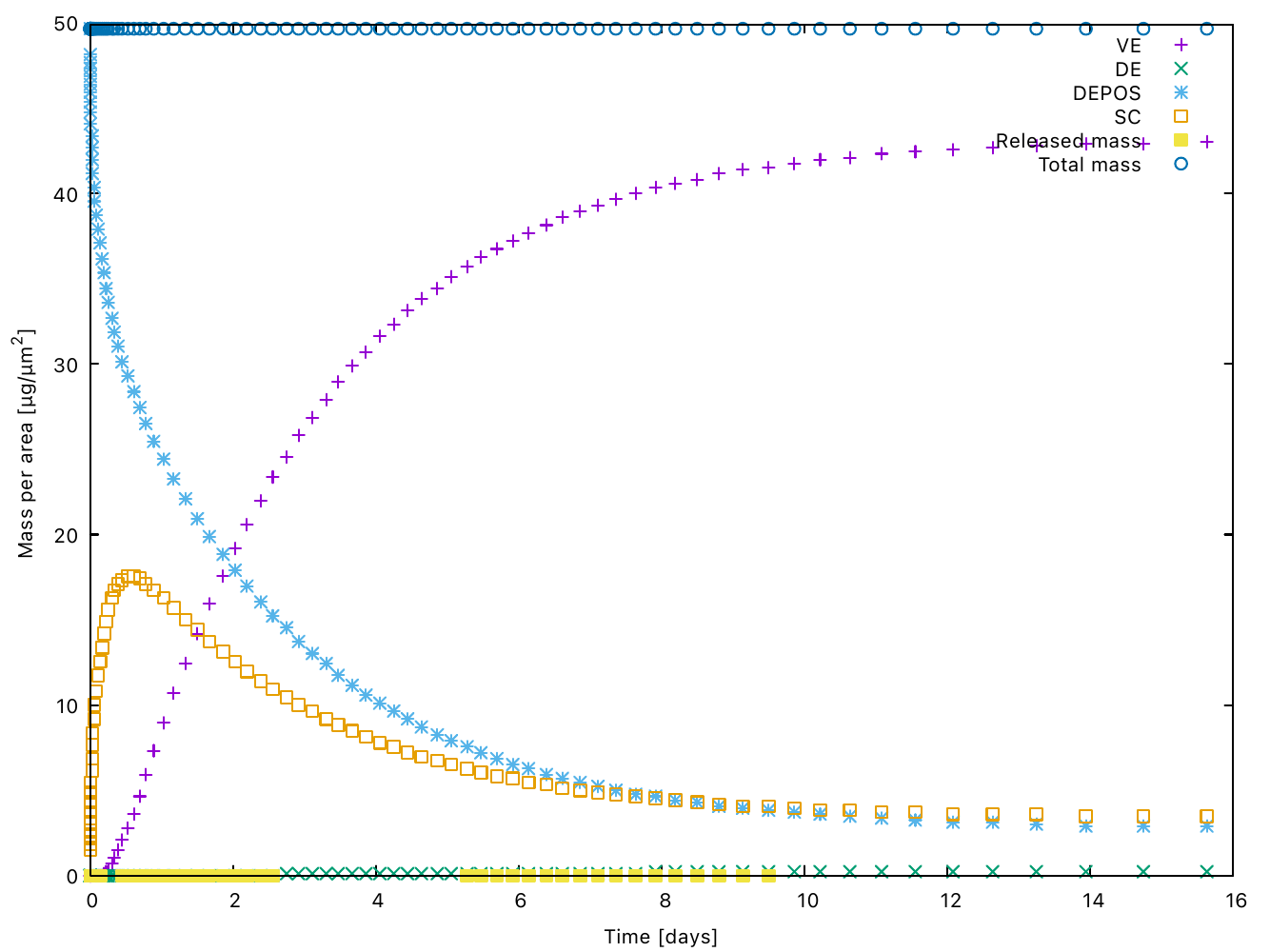}
\caption{Plots of the simulation results for propylparaben (left) and naphthalene (right) in the abdominal skin mesh}  
\label{fig:NP} 
\end{figure}

Naphthalene and propylparaben have similar partition coefficients for SC and VE. However, naphthalene, with a slightly higher partition coefficient ($K_{\mathrm{VE}} = 10.7$) compared to propylparaben ($K_{\mathrm{VE}} = 10.1$), exhibits different absorption behavior due to its distinct diffusion coefficients. Naphthalene has a higher $D$ value than propylparaben, indicating that $D$ significantly influences the absorption rate. Figure \ref{fig:NP} illustrates this effect of $D$ and $K$. While naphthalene has a maximum mass per area of $M_{\mathrm{max}}$ =43.186  \textmu g\slash \textmu m\textsuperscript{2} after 22 days, propylparaben shows a maximum mass per area of $M_{\mathrm{max}}$=18.86 \textmu g\slash \textmu m\textsuperscript{2} after 2 days. Naphthalene exhibits a strong affinity for lipid environments and maintains a low solubility in water, which is indicated by $\log K_{\mathrm{ow}} \approx 3{.}3$. The extended application and exposure time of naphthalene results in significant absorption, as illustrated in the plot. Naphthalene's lipophilicity enables an effective penetration in the SC. Penetration into the deeper tissue, the viable epidermis (VE), is feasible but delayed, as metabolic processes may already be occurring, which is well known for parabens \cite{metab}\cite{jewell}. In contrast, propylparaben is rapidly absorbed in all layers and expelled by the body within 16 days.  .

\subsubsection{Ibuprofen and Resorcinol}
\begin{figure}[htbp]                                    
\centering
\includegraphics[scale=0.22]{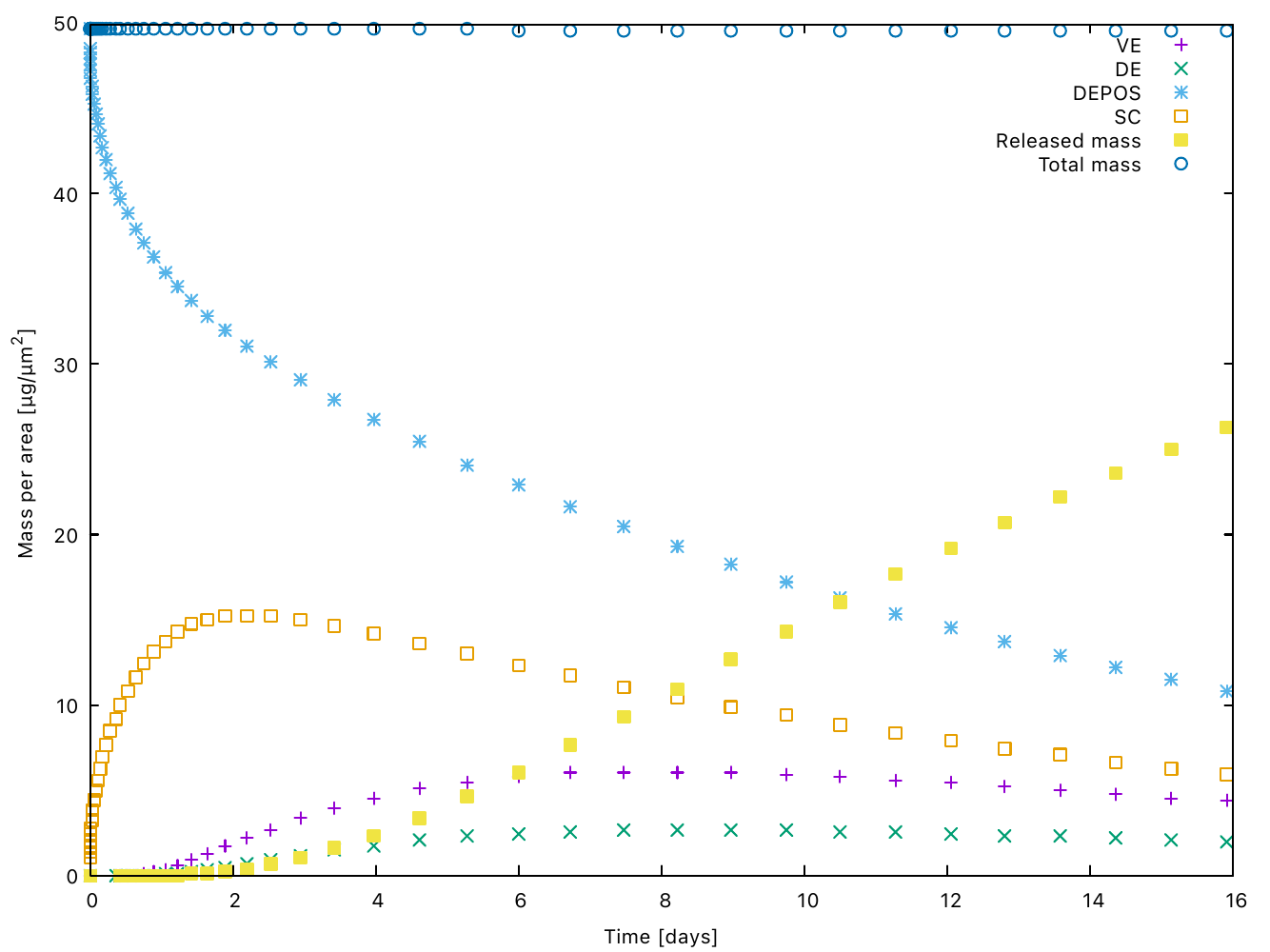}
\includegraphics[scale=0.22]{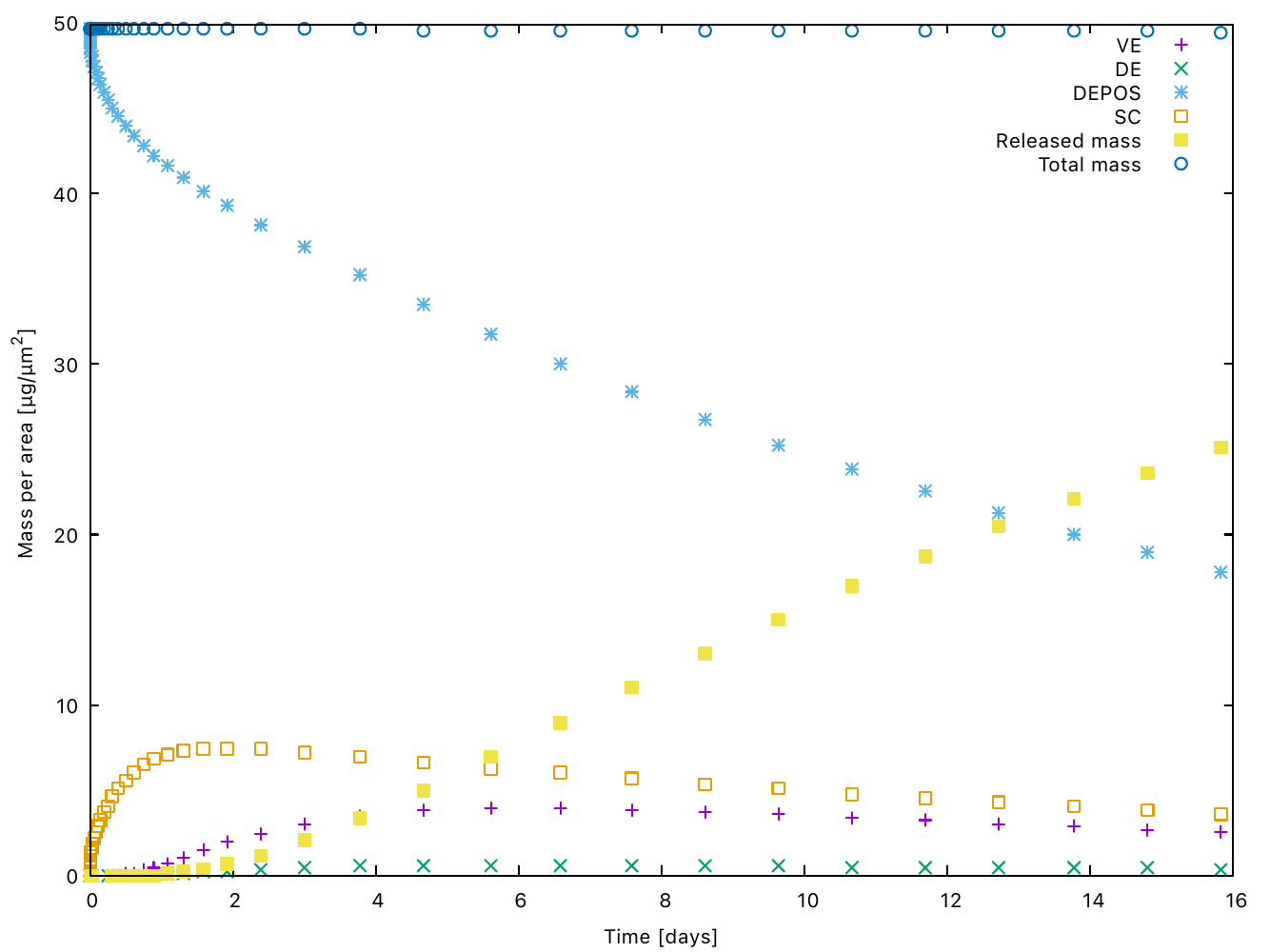}
\caption{Plots of the simulation results for ibuprofen (left) and resorcinol (right) in the abdominal skin mesh}   
\label{fig:IR}
\end{figure}

Topically administered ibuprofen has shown an efficacy comparable to oral ibuprofen in the treatment of joint and soft tissue injuries \cite{migdal}. Ibuprofen exhibits a higher partition coefficient in the VE ($K_{\mathrm{VE}}=6.7$) compared to resorcinol ($K_{\mathrm{VE}}=2.6$). Despite the higher $K_{\mathrm{VE}}$ of ibuprofen the absorption is suboptimal. A nearly equivalent amount of resorcinol is absorbed by the VE. The maximum mass per area for ibuprofen is 5.3 \textmu g\slash \textmu m\textsuperscript{2}, while resorcinol has a maximal mass value of 3.34 \textmu g\slash \textmu m\textsuperscript{2}. In the SC ibuprofen has a higher partition coefficient ($K_{\mathrm{SC}}=7.9$) than resorcinol ($K_{\mathrm{SC}}=3$).  The transfer of resorcinol from the DEPOS layer to the SC occurs at a slower rate (Figure \ref{fig:IR}). Furthermore, based on \cite{anju2024resorcinol}, resorcinol is absorbed through oral, dermal, and subcutaneous routes, followed by rapid metabolism.

\subsubsection{Benzylideneacetone and Benzyl bromide}
\begin{figure}[htbp]   
\begin{center} 
\includegraphics[scale=0.22]{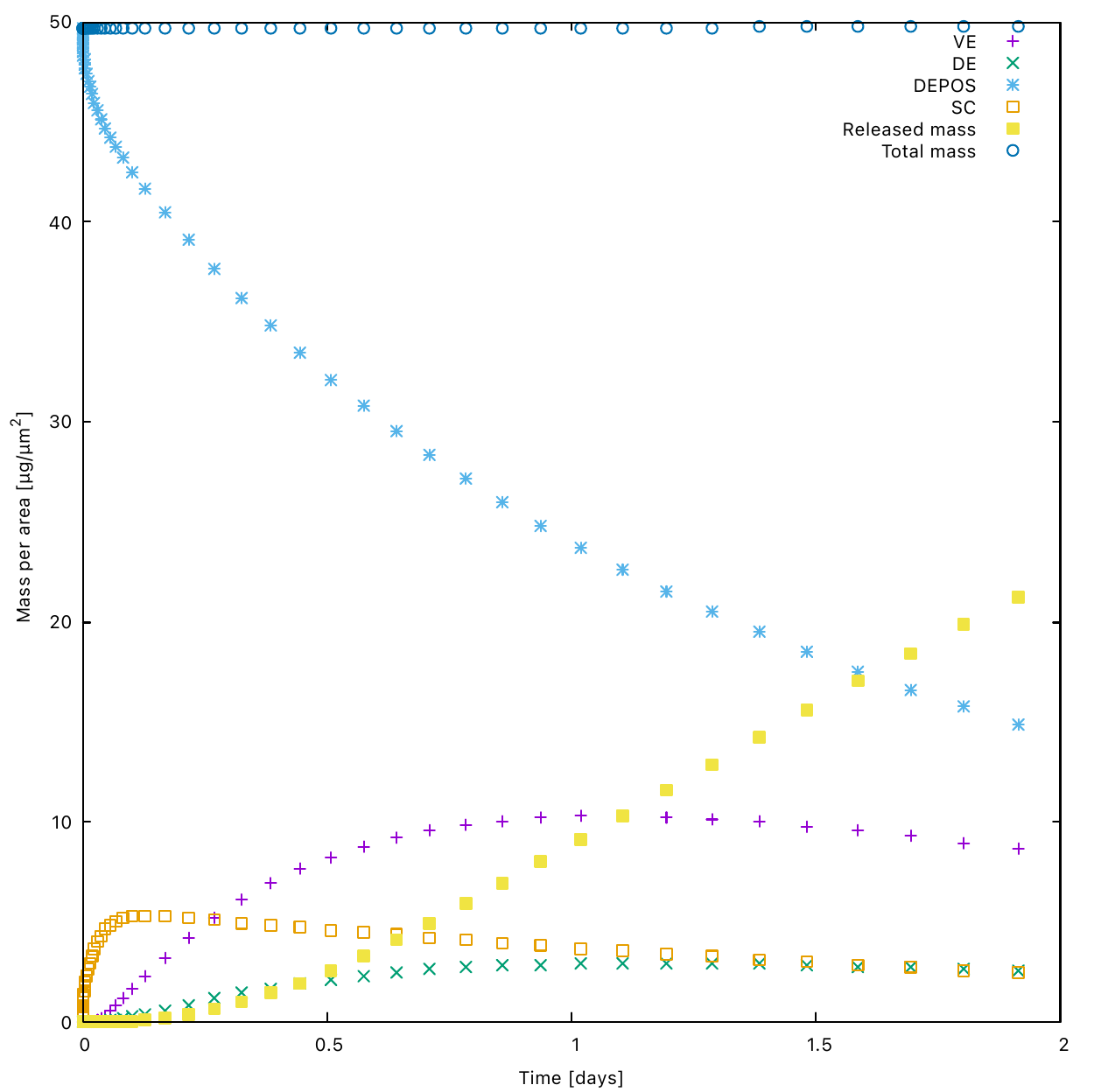}
\includegraphics[scale=0.22]{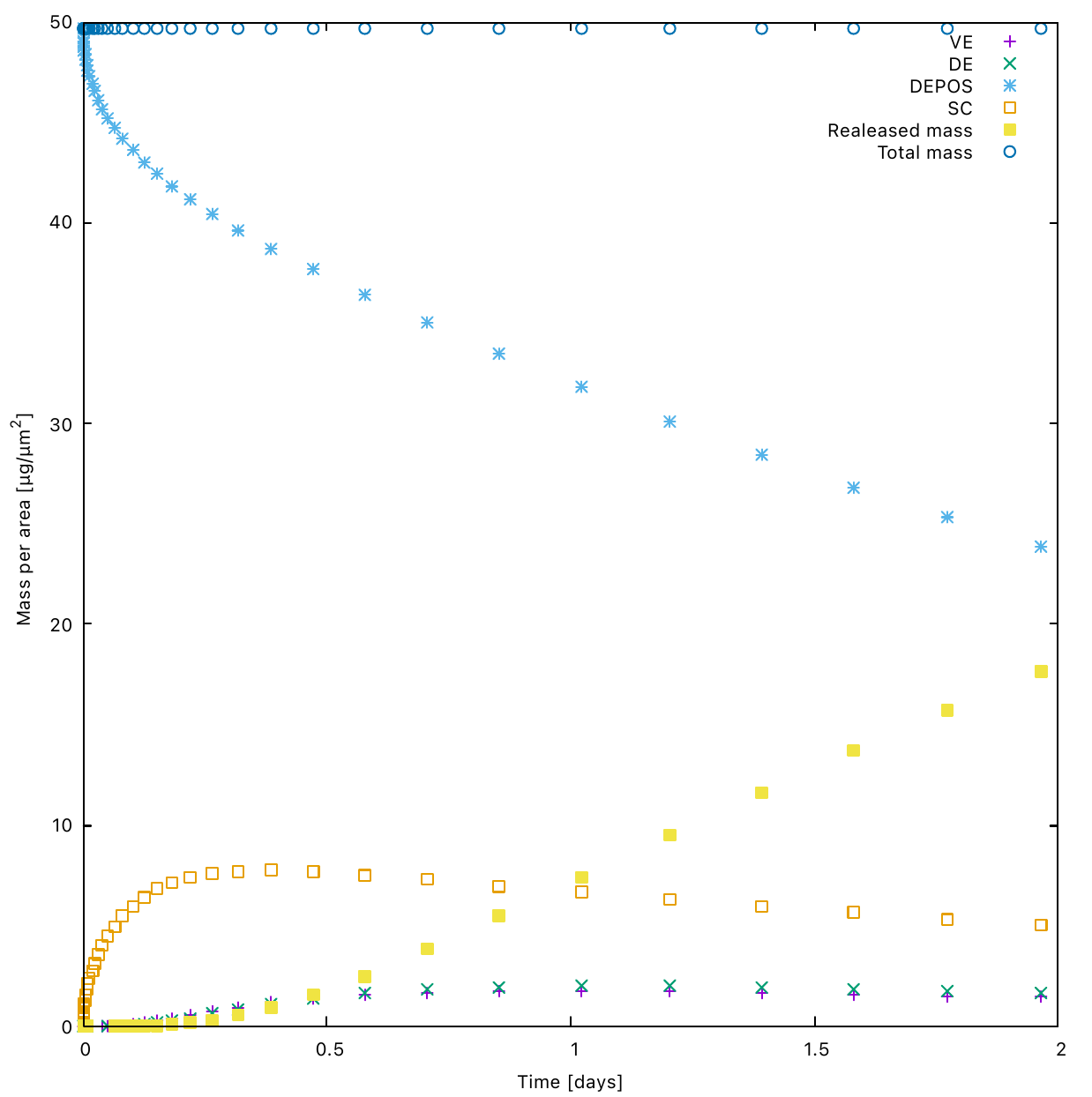}
\caption{Plots of the simulation results for benzylidenacetone (left) and benzyl bromide (right) in the abdominal skin mesh}   
\label{fig:BB}
\end{center}
\end{figure}

In comparison to benzylideneacetone ($K_{\mathrm{SC}} = 0.6$), the SC demonstrates a greater absorption of benzyl bromide,  which can be indicated by the higher K and D ($K_{\mathrm{SC}}$ = 3) (Figure \ref{fig:BB}). The $M_{\mathrm{max}}$ value in the SC for benzyl bromide is approximately $M_{\mathrm{max}} = ~7 $\textmu g\slash \textmu m\textsuperscript{2}, whereas for benzylideneacetone, it is below 5 \textmu g\slash \textmu m\textsuperscript{2}. Overall, the following two observations were made:

\begin{itemize}
\item The maximum mass per area of the SC is greater than that of the VE.
\item The VE exhibits a significantly greater maximum mass per area than the SC.
\end{itemize}

The first observation could be explained by the fact that the SC has a relatively low partition coefficient, causing the chemicals to accumulate in the SC. This accumulation can significantly impact the absorption of the chemical into the VE and DE. The second observation can be explained by the interaction between $D$ and $K$, which can affect the penetration behaviour.  

\subsubsection{p-Chloroaniline and Geraniol}
\begin{figure}[htbp]  
\centering                             
\includegraphics[scale=0.22]{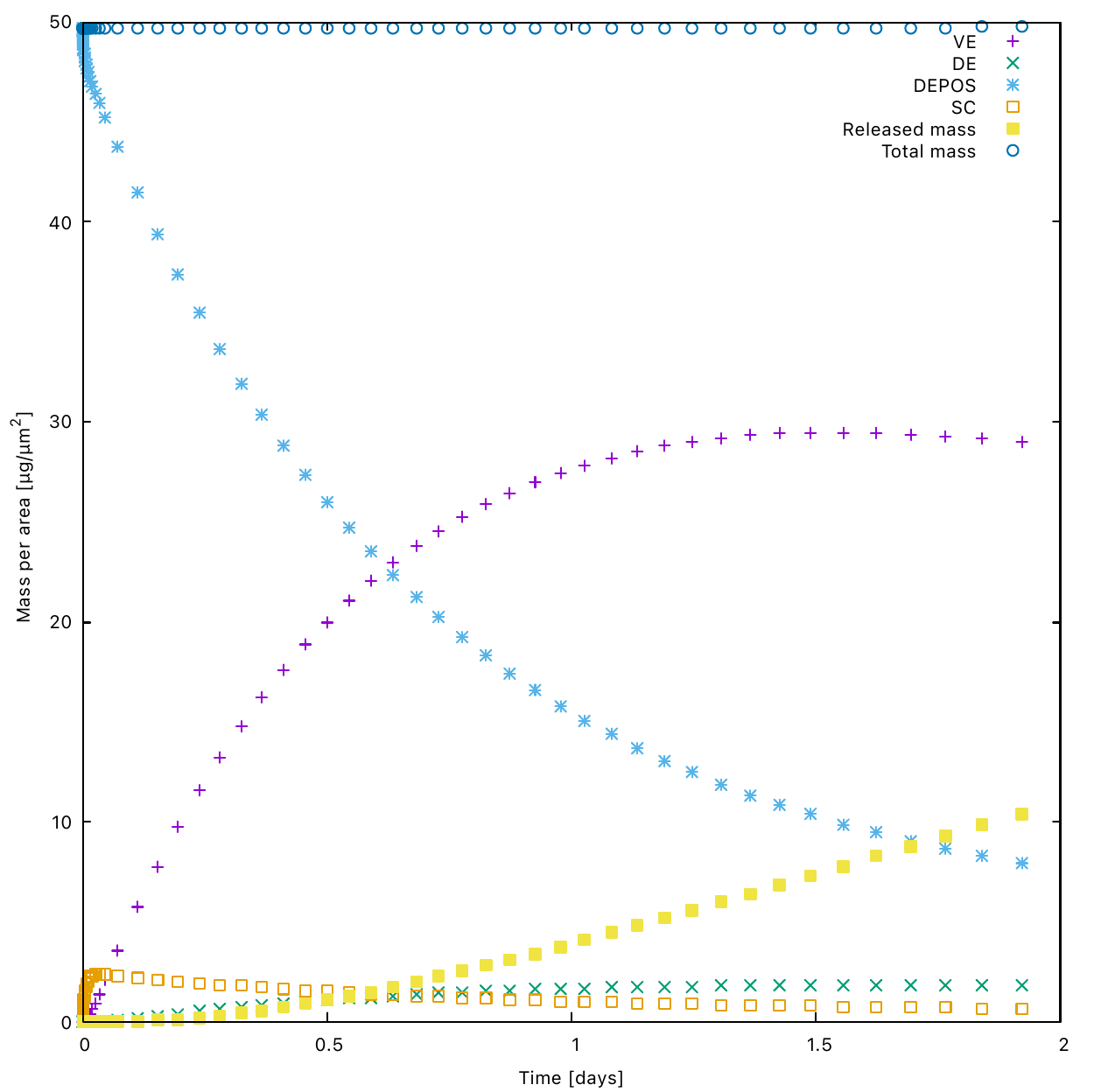}%\\45multi löschen
\includegraphics[scale=0.22]{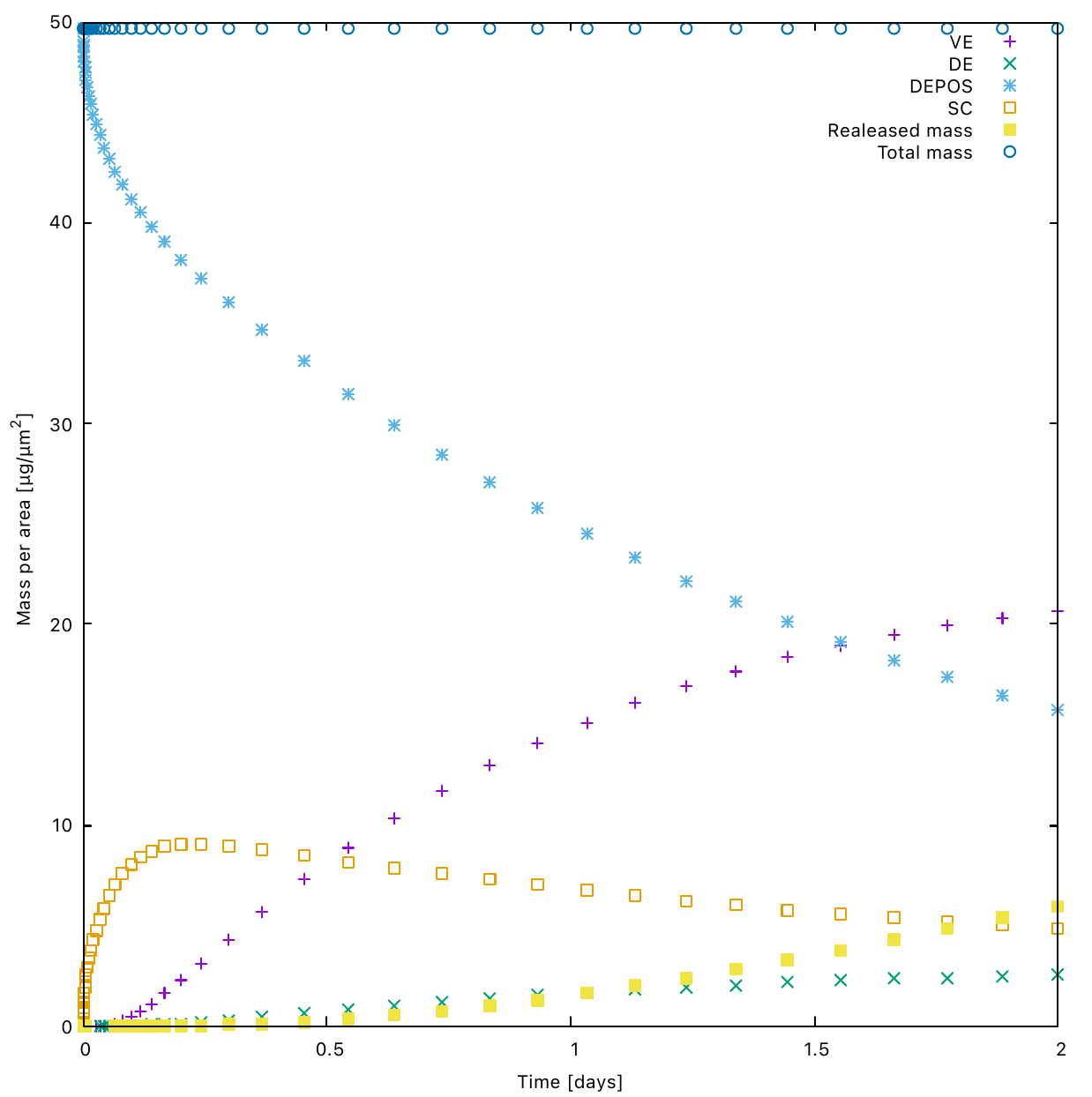}% \\49 multi löschen
\caption{Plots of the simulation results for geraniol (left) and p-chloroaniline (right) in the abdominal skin mesh}   
\label{fig:pG}
\end{figure}

p-Chloroaniline has a moderate lipophilicity, with a $\log K_{\mathrm{ow}}$ value of 1.83 at pH 7.4. This leads to effective permeation into the SC and deeper, hydrophilic layers of the skin, as demonstrated in Figure \ref{fig:pG}. The chemical is absorbed more effectively by the VE than by the SC. In contrast to p-chloroaniline, geraniol is poorly absorbed in the SC. In the VE the $M_{\mathrm{max}}$ value for geraniol is $M_{\mathrm{max}}$ = 29.47 \textmu g\slash \textmu m\textsuperscript{2} after 37 hours. 

\subsubsection{Molecular Weight of Different Chemicals}
\begin{figure}[htbp]                                 
\centering
\includegraphics[scale=0.22]{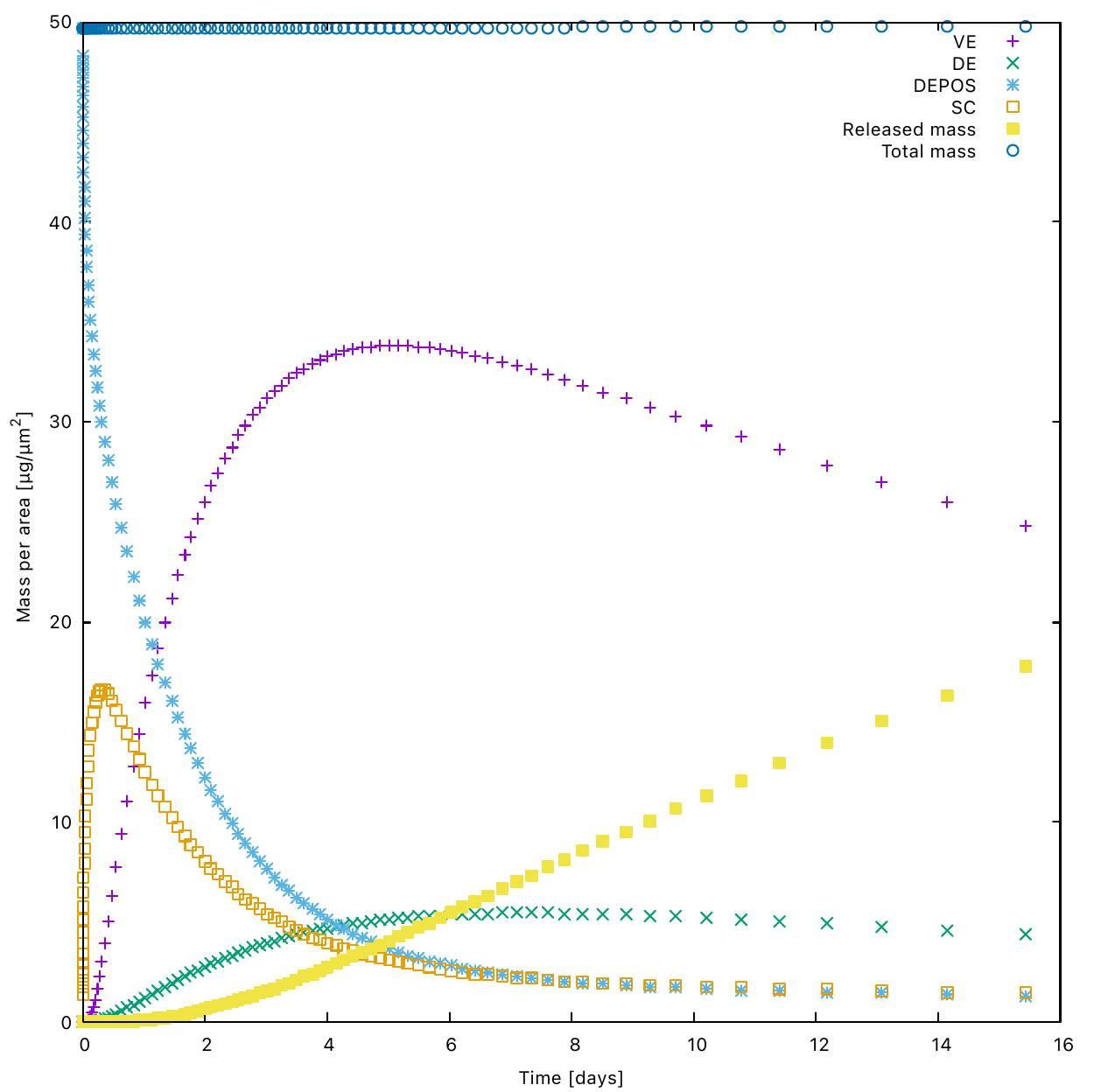}
\includegraphics[scale=0.22]{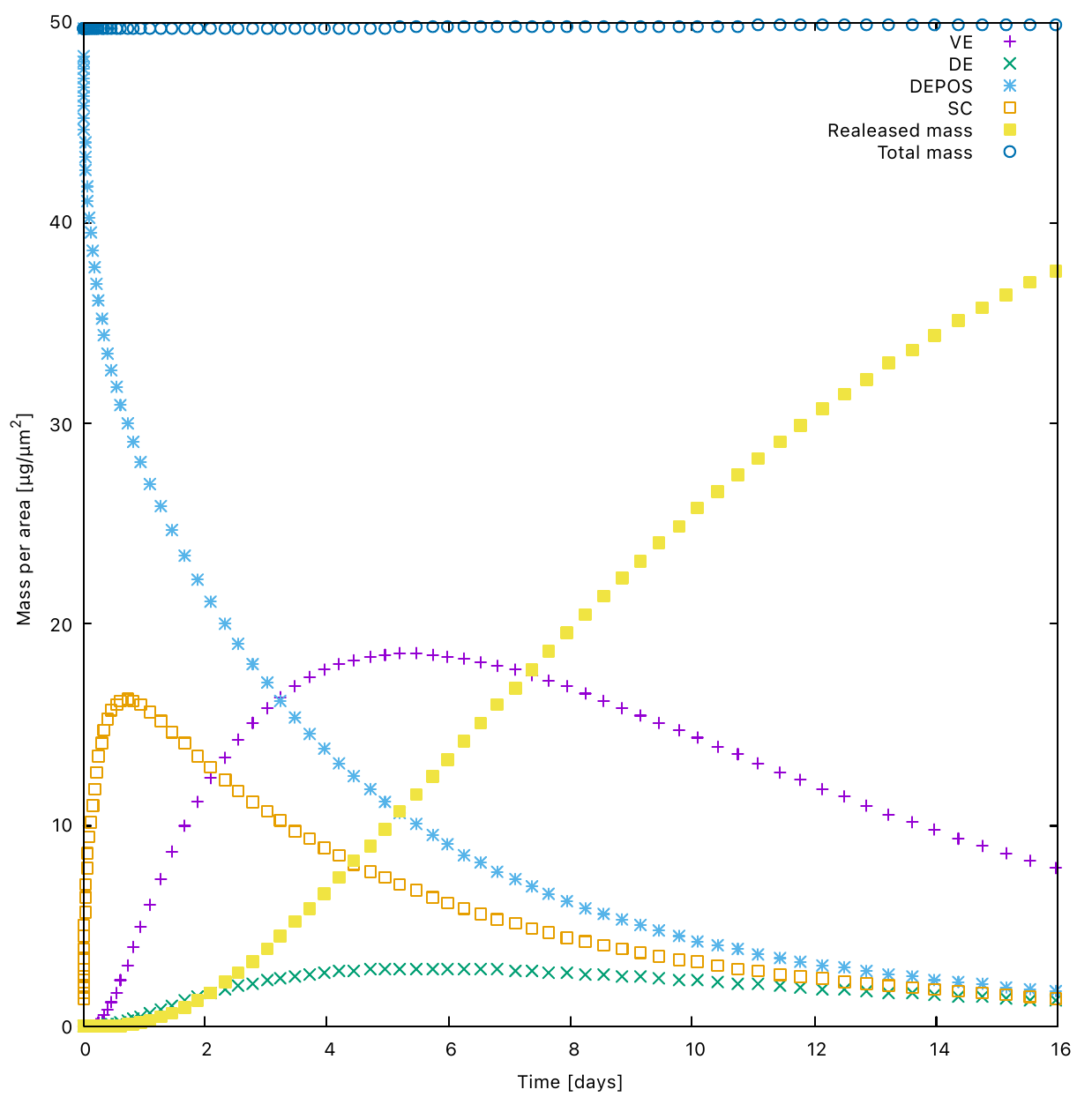}
\caption{Plots of the simulation results for benzophenone (left) and isoeugenol (right) in the abdominal skin mesh}   
\label{fig:NB}
\end{figure}

Figure \ref{fig:NB} shows that benzophenone and isoeugenol are sufficiently absorbed in the SC, VE, and DE. We can observe a faster penetration rate for benzophenone than for isoeugenol. The molecular weight for isoeugenol is 164.201 Da and benzophenone has a higher molecular weight of 182.218 Da. The molecules are small enough to diffuse through the skin immediately, but not too large to have difficulty penetrating the skin. In contrast, some other chemicals lead to different observations. For example, the SC enables the rapid absorption of small molecules, such as benzylideneacetone (MW: 146.186 Da), but exhibits limited diffusivity in VE. The variations in penetration behaviour could be attributed to the distinct properties of the chemical. Chemicals with a molecular weight exceeding 200 Da, including triclosan and basic red 76, have the ability to penetrate the skin effectively.

\subsubsection{Triclosan}
Triclosan is included in many products such as antibacterial soaps, creams, deodorants, or cosmetics \cite{pycke2014human}. Since triclosan has a higher $pK_{\mathrm{a}}(\sim 7.9)$ value than that of the skin pH value ($\sim$4.5--5.5), triclosan can remain longer and persist in skin tissues \cite{dhillon2015triclosan}. As a result of the rapid absorption within the SC, the VE fills up extremely quickly and ultimately transfers its mass to the DE.

\subsubsection{Triclosan in Outerforearm and Chest}
\begin{figure}[htbp]                                 
\centering
\includegraphics[scale=0.22]{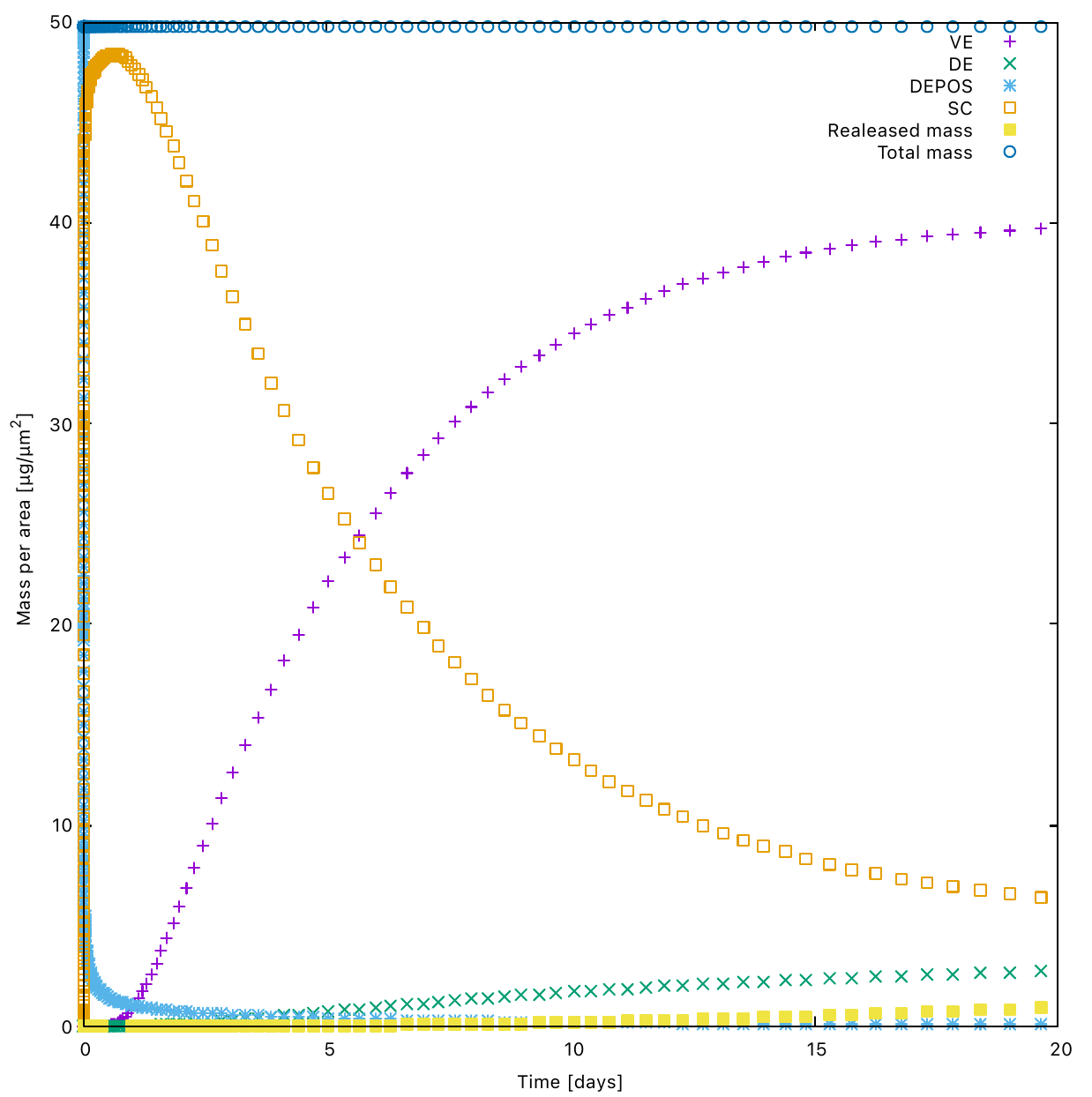}
\includegraphics[scale=0.22]{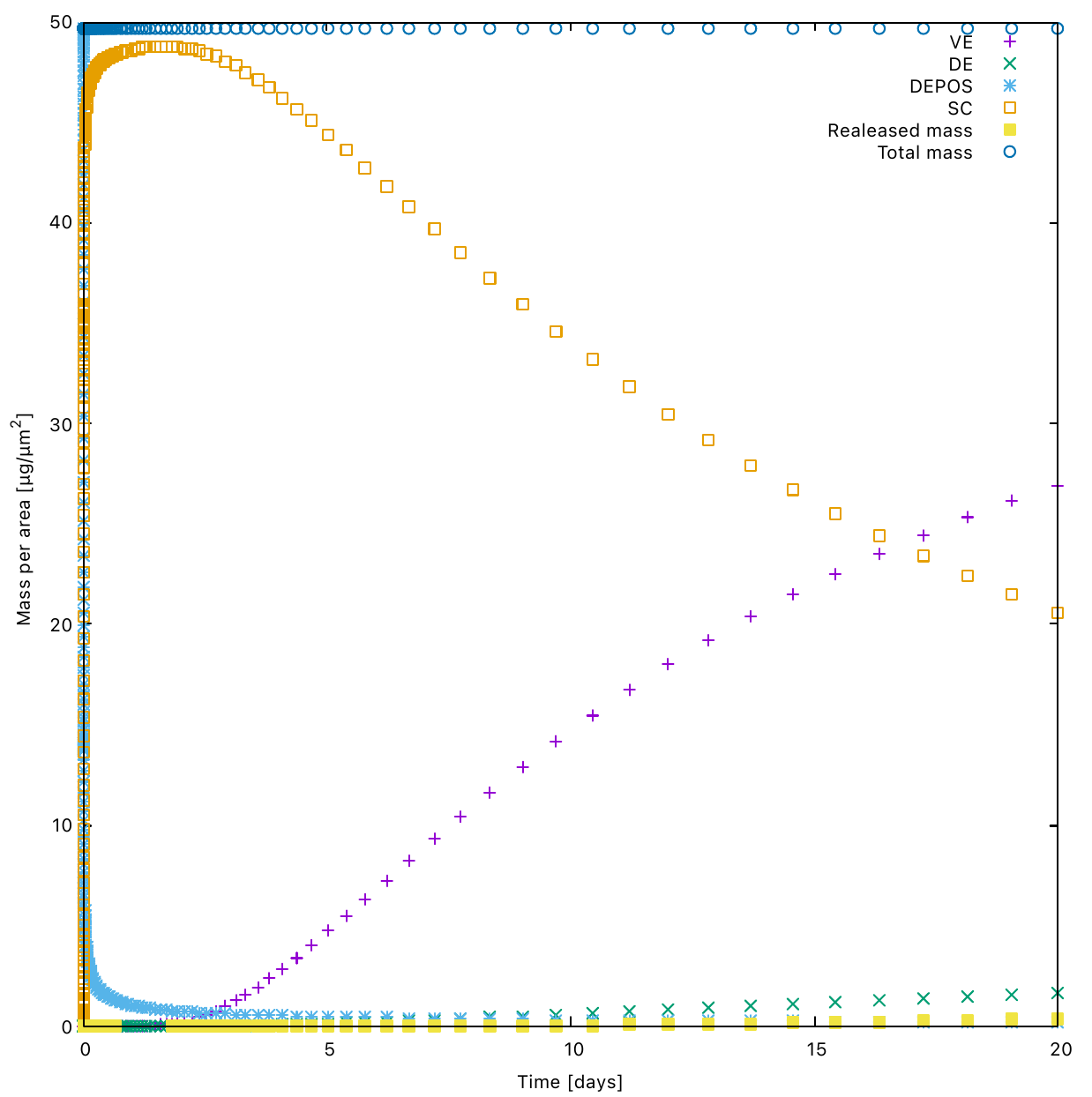}
\caption{Plots showing the simulation results of triclosan in the skin mesh of the chest (left) and the outer forearm (right)}   
\label{fig:TTCO}
\end{figure}

In the beginning, a detailed examination of the different skin thicknesses and regions was discussed. Accordingly, three meshes for distinct skin regions were implemented and used for simulation. Figure \ref{fig:TTCO} presents the simulation results for the skin areas, the chest and the \mbox{outer forearm}.  
\begin{eqnarray}
\tau \approx \frac{L^{2}}{D} \label{zeitab}
\end{eqnarray}
The plots in Figure \ref{fig:TTCO} reveal that the intersection point has shifted to the right by approximately a factor of three. The time estimation for each skin layer is calculated using the formula presented in equation (\ref{zeitab}), where $L$ denotes the thickness of the skin layer and $D$ signifies the diffusion coefficient. The intersection point between SC and VE for the chest region occurs at IP$(5/24.37)$, while for the outer forearm, it is located at IP$(16/23.6)$. 

\subsubsection{Transport of Chemicals through Various Skin Areas}
\begin{figure}[htbp] 
\centering
\subfigure[Chest]{\includegraphics[scale=0.19]{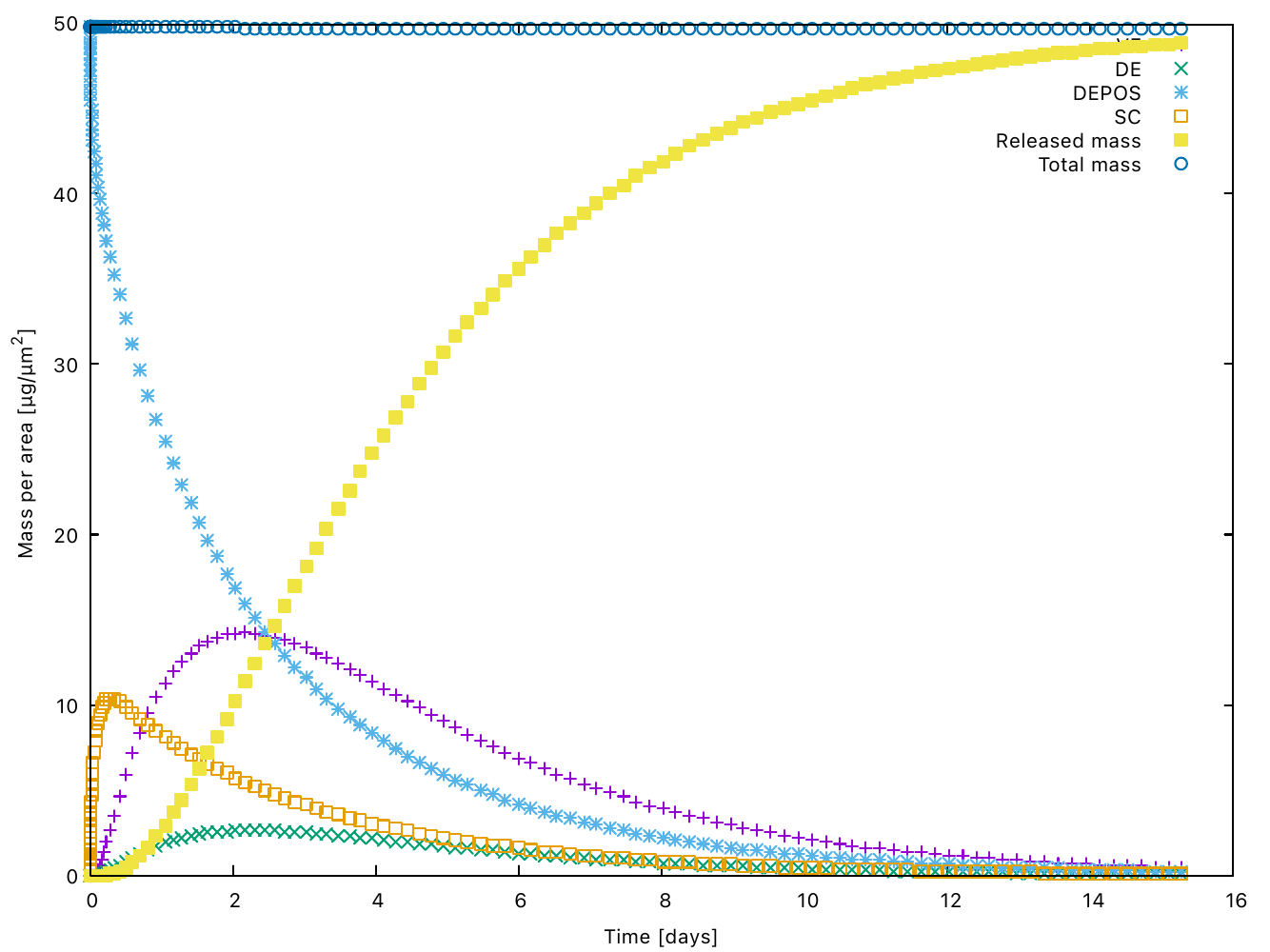}}
\subfigure[Abdomen]{\includegraphics[scale=0.19]{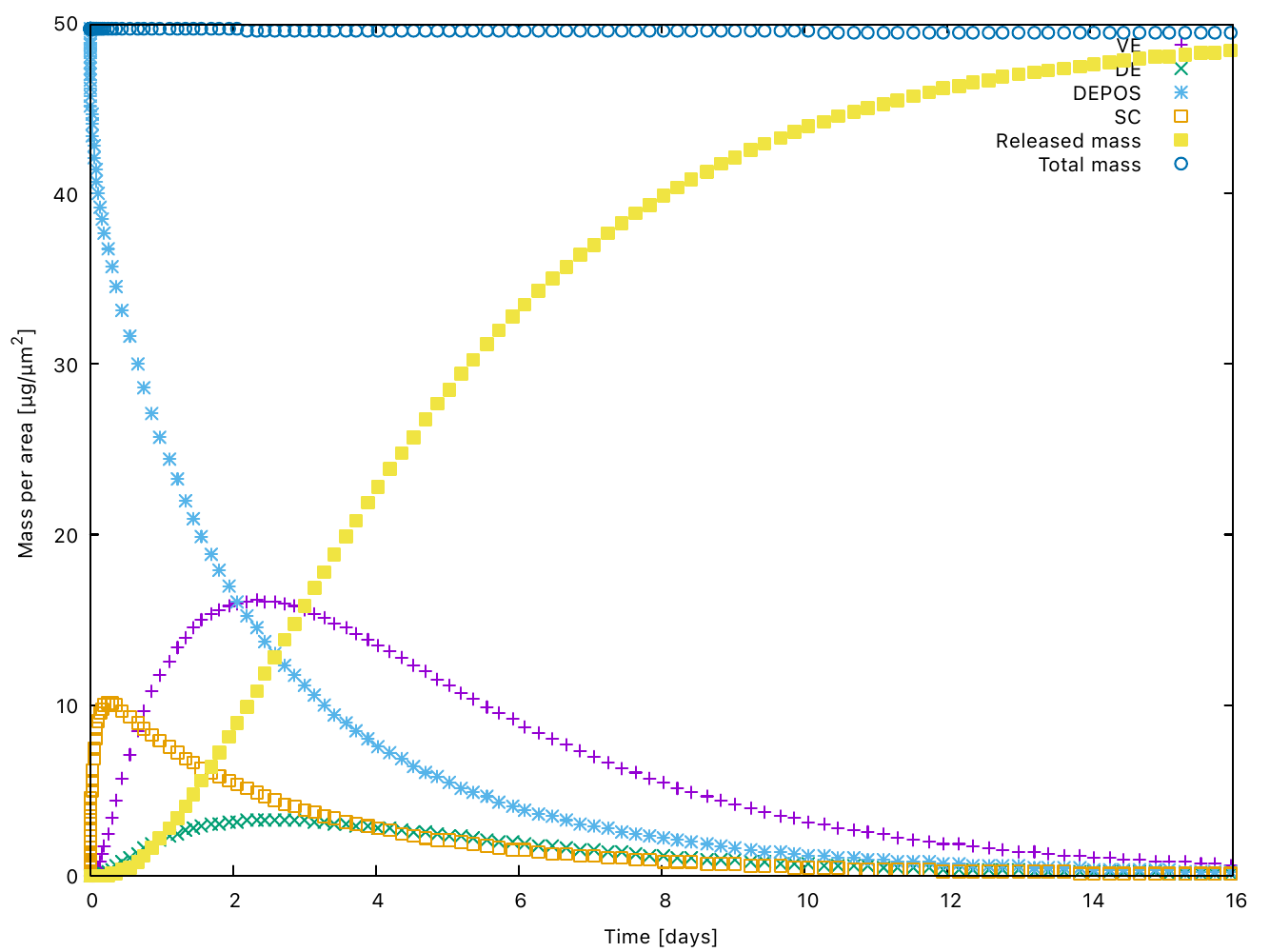}}
\subfigure[Outer forearm]{\includegraphics[scale=0.19]{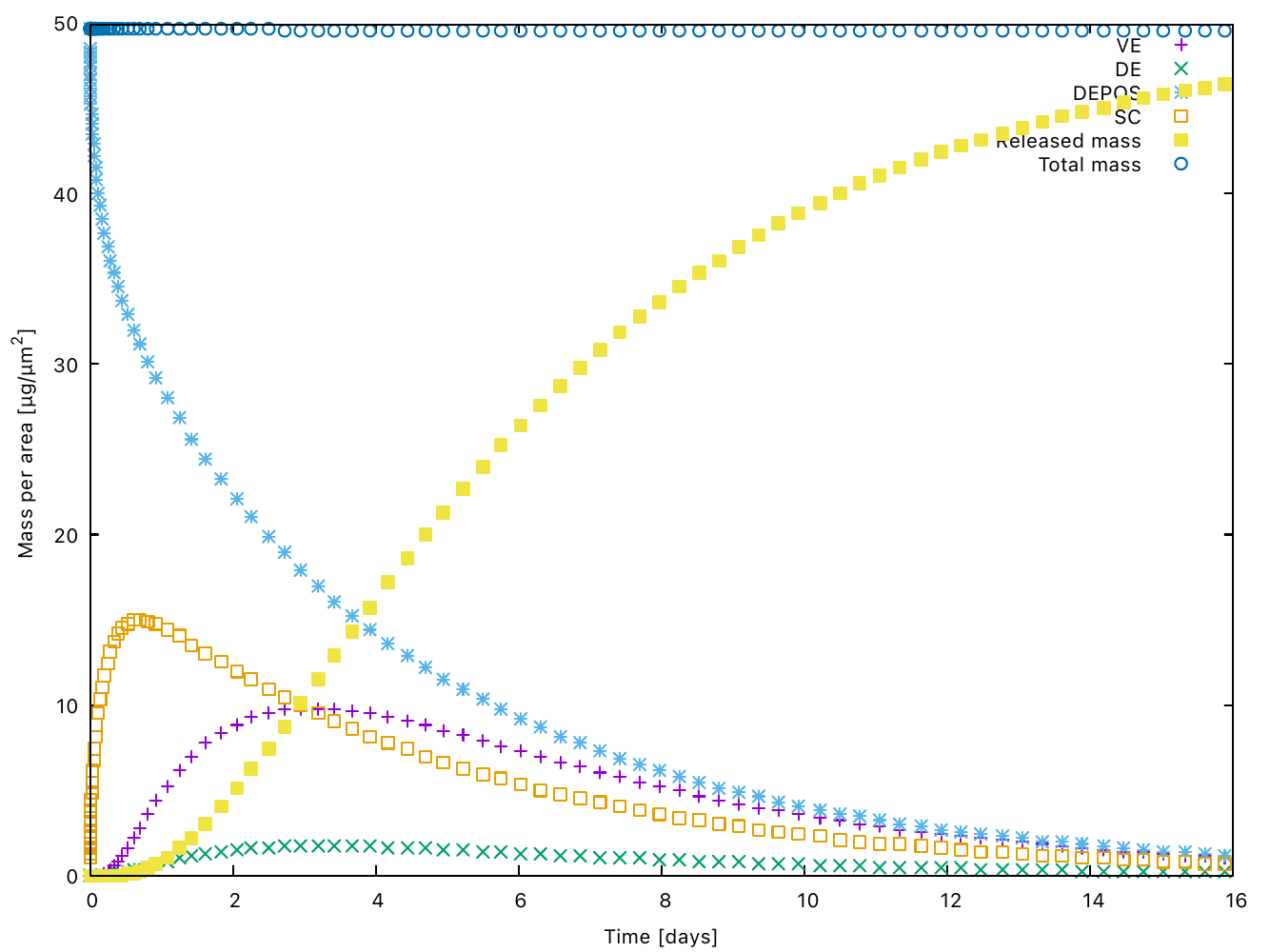}}
\caption{Simulation results of eugenol in a 2D young skin mesh across three distinct skin areas}
\label{fig:grafikco}
\end{figure}

Figure \ref{fig:grafikco} illustrates the penetration of eugenol through the skin regions of the abdomen, chest, and outer forearm. The three plots in Figure \ref{fig:grafikco} show that the SC plot of the outer forearm demonstrates a flatter curve compared to the other two skin regions. The SC of the outer forearm is thicker than that of the chest and abdomen. Therefore, more eugenol can be absorbed. In addition, each of the three skin regions contains distinct intersection points. Regarding the chest region, there are intersection points for the released mass, VE, and the DEPOS layer. In contrast, the abdomen region has an intersection point between the released mass and VE. 

\subsubsection{Simulations}
\begin{figure}[htbp]
\centering   
\includegraphics[scale=0.8]{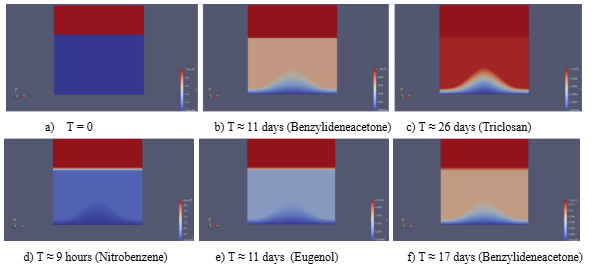}                         
\caption{Simulation results showing the diffusion of chemicals, visualized in ParaView.}   
\label{fig:Para_OLDYOUNG3}
\end{figure}

The simulation results, visualized in ParaView, illustrate the temporal expansion of the chemical distribution across the different skin layers. The concentration profile is represented using a color map, with high concentrations shown in red and low concentrations in blue. In Figure \ref{fig:Para_OLDYOUNG3}a), the DEPOS layer is initially filled at time $t=0$ with a single chemical. The chemical then diffuses into the SC layer and into the deeper layers of the VE and DE. Further simulations in Figure \ref{fig:Para_OLDYOUNG3} illustrate the diffusion of benzylidene acetone, triclosan, eugenol, and nitrobenzene over various time periods. 
\noindent
The quantitative findings of the simulation results are validated by the visual outputs, which illustrate the influence of physicochemical parameters and anatomical variability on skin absorption. These insights can be used to optimize formulation strategies, evaluate exposure risks, and perform experimental validations.
 
\newpage
\section{Conclusion}
This paper aimed to support effective and controlled transdermal permeation of xenobiotics using numerical methods. Therefore, computable skin meshes with realistic anatomical measurements were implemented. By combining structural skin properties with key physicochemical parameters, the simulations provide insights into the pathways and mechanisms of skin permeation. Similar simulation results were observed between 2D and 3D , as well as between young and aged skin meshes. Moreover, the simulated chemicals were significantly influenced by factors such as molecular weight, diffusion, and partition coefficients. These insights might contribute to improving pharmaceutical formulations.\\
In the future, we will enhance the model by integrating capillary loops within the dermal papillae. Integrating them into the simulation will improve the realism of substance transport into the bloodstream, thereby increasing the model's predictive accuracy for systemic absorption.
\newpage

\section*{List of Abbreviations}
\begin{center}
\begin{table}[h!]
\centering
\begin{tabular}{||c || c ||} 
 \hline
 Abbreviation & Authors  \\ [0.5ex] 
 \hline\hline
 C_M & Cinotti_Male   \\ 
 \hline
 C_F & Cinotti_Female   \\
 \hline
 C_T & Cinotti_Total  \\
 \hline
 H & Huzaira   \\
 \hline
 N_O & Neeken_Old   \\
 \hline
 N_Y & Neerken_Young   \\
 \hline
 M_Y & Monnier_Young   \\
 \hline
 SM & Sauermann   \\
 \hline
 S_O & Sandby_Old   \\
 \hline
 S_Y & Sandby_Young   \\ [1ex] 
 \hline
\end{tabular}

\caption{Abbreviations for Figure \ref{fig:TEg}  }
\end{table}
\end{center}
\begin{center}
\begin{table}[h!]
\centering
\begin{tabular}{||c || c ||} 
 \hline
Abbreviation   & Full Form  \\ [0.5ex] 
 \hline\hline
 MW & Molecular weight   \\ 
 \hline
 $D_{DEPOS}$ & Diffusion coefficient of the deposition layer   \\
 \hline
 $D_{SC}$ & Diffusion coefficient of the stratum corneum layer   \\
 \hline
 $D_{VE}$ & Diffusion coefficient of the viable epidermis layer    \\
 \hline
 $D_{DE}$ & Diffusion coefficient of the dermis layer   \\
 \hline
 $K_{DEPOS}$ & Partition coefficient of the deposition layer   \\
 \hline
 $K_{SC}$ & Partition coefficient of the stratum corneum layer   \\
 \hline
 $K_{VE}$ & Partition coefficient of the viable epidermis layer    \\
 \hline
 $K_{DE}$ & Partition coefficient of the dermis layer   \\
 \hline
\end{tabular}
\caption{Abbreviations}
\end{table}
\end{center}
\newpage
\addsec{Acknowledgments}
This research was supported by King Abdullah University of Science and Technology (KAUST).

\addsec{Conflict of interest statement}
The authors declare that they have no known competing financial interests or personal relationships that could have appeared to influence the work reported in this paper.

\addsec{Data availability statement}
All data supporting the findings of this study will be available in the published article and its supplementary materials.

\newpage
\printbibliography

@article{b,
  title={Papillary networks in the dermal--epidermal junction of skin: a biomechanical model},
  author={Ciarletta, Pasquale and Amar, Martine Ben},
  journal={Mechanics Research Communications},
  volume={42},
  pages={68--76},
  year={2012},
}

@article{ga,
  author    = {J. Michael Sorrell and Arnold I. Caplan},
  title     = {Fibroblast heterogeneity: more than skin deep},
  journal   = {Journal of Cell Science},
  year      = {2004},
  volume    = {117},
  number    = {5},
  pages     = {667--675},
}

@article{n7,
  author    = {Monnier, Jilliana and Tognetti, Lorenzo and Miyamoto, Makoto and Suppa, Miriam and Cinotti, Emanuele and Fontaine, Morgane and Perez, Julien and Orte Cano, Carlos and Y{\'e}lamos, Octavio and Puig, Susana and Dubois, Arnaud and Rubegni, Pietro and del Marmol, V{\'e}ronique and Malvehy, Josep and Perrot, Jean-Luc},
  title     = {In vivo characterization of healthy human skin with a novel, non-invasive imaging technique: line--field confocal optical coherence tomography},
  journal   = {Journal of the European Academy of Dermatology and Venereology},
  year      = {2021},
  volume    = {35},
  number    = {9},
  pages     = {1966--1974},
}

@article{n4,
  author    = {Sauermann, Kirsten and Clemann, Sven and Jaspers, S{\"o}ren and Gambichler, Thilo and Altmeyer, Peter and Hoffmann, Klaus and Ennen, Joachim},
  title     = {Age-related changes of human skin investigated with histometric measurements by confocal laser scanning microscopy {in vivo}},
  journal   = {Skin Research and Technology},
  year      = {2002},
  volume    = {8},
  number    = {1},
  pages     = {52--56},
}

@article{n77,
  author    = {Lee, Y. and Hwang, K.},
  title     = {Skin thickness of Korean adults},
  journal   = {Surgical and Radiologic Anatomy},
  year      = {2002},
  volume    = {24},
  number    = {3--4},
  pages     = {183--189},
}

@article{g,
  author    = {Girardeau, Sarah and Mine, Sol{\`e}ne and Pageon, Herv{\'e} and Asselineau, Daniel},
  title     = {The Caucasian and African skin types differ morphologically and functionally in their dermal component},
  journal   = {Experimental Dermatology},
  year      = {2009},
  volume    = {18},
  number    = {8},
  pages     = {704--711},
}

@article{n3,
  author    = {Huzaira, Misbah and Rius, Francisca and Rajadhyaksha, Milind and Anderson, R. Rox and Gonz{\'a}lez, Salvador},
  title     = {Topographic variations in normal skin, as viewed by {in vivo} reflectance confocal microscopy},
  journal   = {Journal of Investigative Dermatology},
  year      = {2001},
  volume    = {116},
  number    = {6},
  pages     = {846--852},
}

@online{ug4,
  title     = {UG4},
  url       = {http://ug4.github.io/docs/page_u_g4_introduction.html},
  urldate   = {2025-03-09},
}

@techreport{ux,
  author      = {Arne N{\"a}gel and Peter Deuflhard and Gabriel Wittum},
  title       = {Efficient Stiff Integration of Density Driven Flow Problems},
  institution = {ZIB},
  address     = {Takustr. 7, 14195 Berlin},
  number      = {18-54},
  language    = {eng},
  year        = {2018}
}

@book{uu,
  author={Peter Deuflhard and Martin Weiser},
  title={Adaptive Numerical Solutions of PDEs}, 
  year={2012} ,
  publisher={De Gruyter},
  adress={Berlin},
}

@article{Tfentanyl,
  title={Predicting transdermal fentanyl delivery using mechanistic simulations for tailored therapy},
  author={Defraeye, Thijs and Bahrami, Flora and Ding, Lu and Malini, Riccardo Innocenti and Terrier, Alexandre and Rossi, Ren{\'e} M},
  journal={Frontiers in pharmacology},
  volume={11},
  pages={585393},
  year={2020},
}

@article{n6,
  title={Effect of 3D microstructure of dermal papillae on SED concentration at a mechanoreceptor location},
  author={Pham, Trung Quang and Hoshi, Takayuki and Tanaka, Yoshihiro and Sano, Akihito},
  journal={PloS one},
  volume={12},
  number={12},
  year={2017},
}

@book{mdiff,
  title={Diffusion: mass transfer in fluid systems},
  author={Cussler, Edward Lansing},
  year={2009},
  publisher={Cambridge university press}
}

@book{hackbusch,
  title={Multi-grid methods and applications},
  author={Hackbusch, Wolfgang},
  volume={4},
  year={2013},
  publisher={Springer Science \& Business Media}
}

@online{proM,
  author       = {Reiter, Sebastian},
  title        = {{ProMesh}},
  year         = {2022},
  URL = {http://www.promesh3d.com},
  URLDATE = {2025-07-04},
}

@article{ellison,
  title={Partition coefficient and diffusion coefficient determinations of 50 compounds in human intact skin, isolated skin layers and isolated stratum corneum lipids},
  author={Ellison, Corie A and Tankersley, Kevin O and Obringer, Cindy M and Carr, Greg J and Manwaring, John and Rothe, Helga and Duplan, H{\'e}l{\`e}ne and G{\'e}ni{\`e}s, Camille and Gr{\'e}goire, S{\'e}bastien and Hewitt, Nicola J and others},
  journal={Toxicology in Vitro},
  volume={69},
  pages={104990},
  year={2020},
}

@article{dancik,
  title={Design and performance of a spreadsheet-based model for estimating bioavailability of chemicals from dermal exposure},
  author={Dancik, Yuri and Miller, Matthew A and Jaworska, Joanna and Kasting, Gerald B},
  journal={Advanced drug delivery reviews},
  volume={65},
  number={2},
  pages={221--236},
  year={2013},
}

@article{selzer,
  title={Finite dose skin mass balance including the lateral part: comparison between experiment, pharmacokinetic modeling and diffusion models},
  author={Selzer, D and Hahn, T and N\"agel, A and Heisig, M and Kostka, KH and Lehr, CM and Neumann, D and Schaefer, UF and Wittum, Gabriel},
  journal={Journal of Controlled Release},
  volume={165},
  number={2},
  pages={119--128},
  year={2013},
}

@article{Thio,
  title={Final amended report on the safety assessment of ammonium thioglycolate, butyl thioglycolate, calcium thioglycolate, ethanolamine thioglycolate, ethyl thioglycolate, glyceryl thioglycolate, isooctyl thioglycolate, isopropyl thioglycolate, magnesium thioglycolate, methyl thioglycolate, potassium thioglycolate, sodium thioglycolate, and thioglycolic acid},
  author={Burnett, Christina L and Bergfeld, Wilma F and Belsito, Donald V and Klaassen, Curtis D and Marks, James G and Shank, Ronald C and Slaga, Thomas J and Snyder, Paul W and Andersen, F Alan},
  journal={International journal of toxicology},
  volume={28},
  number={4\_suppl},
  pages={68--133},
  year={2009},
}

@article{triclosan,
  title={Triclosan exposure, transformation, and human health effects},
  author={Weatherly, Lisa M and Gosse, Julie A},
  journal={Journal of Toxicology and Environmental Health, Part B},
  volume={20},
  number={8},
  pages={447--469},
  year={2017},
}

@phdthesis{nitro,
  title={Kosmetische Mittel vom Kaiserreich bis zur Zeit der Weimarer Republik: Herstellung, Entwicklung und Verbraucherschutz},
  author={Riewerts, Kerrin},
  year={2005},
  school={Staats-und Universit{\"a}tsbibliothek Hamburg Carl von Ossietzky}
}

@article{acetone,
  title={Inhibition effects of benzylideneacetone, benzylacetone, and 4-phenyl-2-butanol on the activity of mushroom tyrosinase},
  author={Liu, Xuan and Jia, Yu-long and Chen, Jing-wei and Liang, Ge and Guo, Hua-yun and Hu, Yong-hua and Shi, Yan and Zhou, Han-Tao and Chen, Qing-Xi},
  journal={Journal of bioscience and bioengineering},
  volume={119},
  number={3},
  pages={275--279},
  year={2015},
}

@article{cherian,
  title={Safety Assessment of Basic Red 76 as Used in Cosmetics},
  author={Cherian, Priya and Bergfeld, Wilma F and Belsito, Donald V and Hill, Ronald A and Klaassen, Curtis D and Liebler, Daniel C and Marks Jr, James G and Shank, Ronald C and Slaga, Thomas J and Snyder, Paul W and others},
  journal={International journal of toxicology},
  volume={43},
  number={3\_suppl},
  pages={109S--119S},
  year={2024},
}

@article{propyl,
  title={Toxic effects of paraben and its relevance in cosmetics: a review},
  author={Alaba, Pierce Adrianne A and Canete, Edjay D and Pantalan, Bin Salih S and Taguba, Joanna Marie C and Yu, Lovely Dianne I and Faller, Erwin M},
  journal={Int J Res},
  volume={2582},
  pages={7421},
  year={2022}
}

@article{propyl2,
  title={Safety assessment of propyl paraben: a review of the published literature},
  author={Soni, MG and Burdock, GA and Taylor, SL and Greenberg, NA},
  journal={Food and Chemical Toxicology},
  volume={39},
  number={6},
  pages={513--532},
  year={2001},
}

@article{Ethy2,
  author={Elmets, Craig A and Yusuf, Nabiha},
  title={Murine skin carcinogenesis and the role of immune system dysregulation in the tumorigenicity of 2-ethylhexyl acrylate},
  journal={Biomedicine hub},
  volume={5},
  number={2},
  pages={1--16},
  year={2020},
}

@article{Ethy21,
  author    = {R. P. Wenzel-Hartung and H. Brune and H.-J. Klimisch},
  title     = {Dermal oncogenicity study of 2-ethylhexyl acrylate by epicutaneous application in male C3H/HeJ mice},
  journal   = {Journal of Cancer Research and Clinical Oncology},
  volume    = {115},
  number    = {6},
  pages     = {543--549},
  year      = {1989}
}

@online{Naph0,
  author       = {Gervais, J. and Luukinen, B. and Buhl, K. and Stone, D.},
  title        = {Naphthalene General Fact Sheet},
  year         = {2010},
  organization = {National Pesticide Information Center, Oregon State University Extension Services},
  url          = {http://npic.orst.edu/factsheets/naphgen.html},
  urldate      = {2025-09-06}
}

@article{n1,
  title={Structural skin changes in elderly people investigated by reflectance confocal microscopy},
  author={Cinotti, E and Bovi, C and Tonini, G and Labeille, B and Heus{\`e}le, C and Nizard, C and Schnebert, S and Aubailly, S and Barth{\'e}l{\'e}my, JC and Cambazard, F and others},
  journal={Journal of the European Academy of Dermatology and Venereology},
  volume={34},
  number={11},
  pages={2652--2658},
  year={2020},
}

@online{Rc,
  author       = {{DrugBank Online}},
  title        = {Resorcinol},
  year         = {2015},
  URL = {https://go.drugbank.com/drugs/DB11085},
  URLDATE = {2025-10-07},
}

@article{Neil,
  title={Human skin explant model for the investigation of topical therapeutics},
  author={Neil, Jessica E and Brown, Marc B and Williams, Adrian C},
  journal={Scientific Reports},
  volume={10},
  number={1},
  pages={21192},
  year={2020},
}

@article{blair2020skin,
  title={Skin structure--function relationships and the wound healing response to intrinsic aging},
  author={Blair, Michael J and Jones, Jake D and Woessner, Alan E and Quinn, Kyle P},
  journal={Advances in wound care},
  volume={9},
  number={3},
  pages={127--143},
  year={2020},
}

@article{OY,
  title={Characterization of age-related effects in human skin: a comparative study that applies confocal laser scanning microscopy and optical coherence tomography},
  author={Neerken, Sieglinde and Lucassen, Gerald W and Bisschop, Marielle A and Lenderink, Egbert and Nuijs, Tom},
  journal={Journal of biomedical optics},
  volume={9},
  number={2},
  pages={274--281},
  year={2004},
}

@article{vybohova2012quantitative,
  title={Quantitative changes of the capillary bed in aging human skin},
  author={V{\`y}bohov{\'a}, Desanka and Mellov{\'a}, Yvetta and Adamicov{\'a}, Katar{\'\i}na and Adamkov, Mari{\'a}n and He{\v{s}}kov{\'a}, Gabriela},
  journal={Histology and Histopathology},
  volume={27},
  number={7},
  pages={961},
  year={2012}
}

@techreport{jewell,
  title        = {Scientific Committee on Consumer Safety (SCCS) Opinion on Propylparaben (PP)},
  author       = {{Scientific Committee on Consumer Safety (SCCS)}},
  year         = {2021},
  institution  = {European Commission},
  number       = {SCCS/1623/20},
  url          = {https://health.ec.europa.eu/system/files/2022-08/sccs_o_243.pdf}
}

@article{migdal,
  title={Topical administration of ibuprofen for injured athletes: considerations, formulations, and comparison to oral delivery},
  author={Manoukian, Martin Anthony Christopher and Migdal, Christopher William and Tembhekar, Amode Ravindra and Harris, Jerad Alexander and DeMesa, Charles},
  journal={Sports medicine-open},
  volume={3},
  number={1},
  pages={36},
  year={2017},
}

@article{metab,
  title     = {Exposure to naphthalene induces naphthyl-keratin adducts in human epidermis in vitro and in vivo},
  author    = {Kang-Sickel, Juei-Chuan C. and Stober, Vandy P. and French, John E. and Nylander-French, Leena A.},
  journal   = {Biomarkers},
  volume    = {15},
  number    = {6},
  pages     = {488--497},
  year      = {2010},
}

@article{kretsos,
  title={A geometrical model of dermal capillary clearance},
  author={Kretsos, Kosmas and Kasting, Gerald B},
  journal={Mathematical biosciences},
  volume={208},
  number={2},
  pages={430--453},
  year={2007},
  publisher={Elsevier}
}

@article{baqar2019dermal,
  title={Microcirculation disorders of the skin},
  author={Creutzig, A and Caspary, L},
  journal={Der Internist},
  volume={35},
  number={6},
  pages={546--556},
  year={1994}
}

@article{anju2024resorcinol,
  title={Resorcinol in Dermatology},
  author={Anju George, C and Chhabra, Namrata},
  journal={Indian Journal of Postgraduate Dermatology Volume},
  volume={2},
  number={1},
  pages={20},
  year={2024}
}

@article{pycke2014human,
  title     = {Human fetal exposure to triclosan and triclocarban in an urban population from Brooklyn, New York},
  author    = {Pycke, Benny F. G. and Geer, Laura A. and Dalloul, Mudar and Abulafia, Ovadia and Jenck, Alizee M. and Halden, Rolf U.},
  journal   = {Environmental Science \& Technology},
  volume    = {48},
  number    = {15},
  pages     = {8831--8838},
  year      = {2014},
}

@article{dhillon2015triclosan,
  title={Triclosan: current status, occurrence, environmental risks and bioaccumulation potential},
  author={Dhillon, Gurpreet Singh and Kaur, Surinder and Pulicharla, Rama and Brar, Satinder Kaur and Cled{\'o}n, Maximiliano and Verma, Mausam and Surampalli, Rao Y},
  journal={International journal of environmental research and public health},
  volume={12},
  number={5},
  pages={5657--5684},
  year={2015},
}

@inbook{yousef2017anatomy,
  author    = {Hani Yousef and Mandy Alhajj and Sandeep Sharma},
  title     = {Anatomy, Skin (Integument), Epidermis},
  booktitle = {StatPearls},
  year      = {2017},
  publisher = {StatPearls Publishing},
}

@book{knabner2013numerik,
  title={Numerik partieller Differentialgleichungen: eine anwendungsorientierte Einf{\"u}hrung},
  author={Knabner, Peter and Angermann, Lutz},
  year={2013},
  publisher={Springer-Verlag}
}

@article{kruczek2005analytical,
  title={Analytical solution for the effective time lag of a membrane in a permeate tube collector in which Knudsen flow regime exists},
  author={Kruczek, B and Frisch, HL and Chapanian, R},
  journal={Journal of membrane Science},
  volume={256},
  number={1-2},
  pages={57--63},
  year={2005},
}

\section*{Figures }
\label{app1}
\label{sec:figures}
 
\newpage
\subsection*{Plots of Simulation Results}
\label{subsec:Simresults}

%%%%%%%%%%%%%%%%%%%%%%%%%%%%%%%%%%%%%%%%%%%%%%%%%%%%%%%%%%%%

%\addxcontentsline{toc}{section}{List of Figures}
%\enlargethispage{\baselineskip}
%\listoffigures % Abbildungen
%\addxcontentsline{toc}{section}{List of Tables}
%\listoftables % Tabellen 

\end{document}